\documentclass[12pt]{article}
\usepackage[utf8]{inputenc}
\usepackage{natbib}
\RequirePackage{amsthm,amsmath,amsfonts,amssymb}
\usepackage[margin=1in]{geometry}
\usepackage{array,epsfig,rotating}
\usepackage{color}
\usepackage[dvipsnames]{xcolor}
\usepackage{hyperref}
\hypersetup{
  colorlinks,
  linkcolor={red!50!black},
  citecolor={blue},
  urlcolor={blue!80!black}
}

\newtheorem{assumption}{Assumption}

\newtheorem{definition}{Definition}
\newtheorem{theorem}{Theorem}
\newtheorem{proposition}{Proposition}
\newtheorem{corollary}{Corollary}
\newtheorem{lemma}{Lemma}

\usepackage{tikz}
\usetikzlibrary{shapes, calc, shapes, arrows, positioning}
\tikzstyle{qedge}=[->,thick,black]
\definecolor{myblue}{rgb}{0.0265,    0.6137,    0.8135}
\definecolor{myyellow}{rgb}{0.9290,    0.6940,    0.1250}
\tikzstyle{neuron}=[draw,circle,minimum size=26pt,inner sep=0pt, fill=black!10]
\tikzstyle{hidden}=[draw,circle,minimum size=26pt,inner sep=0pt, fill=white]
\usetikzlibrary{arrows.meta}
\tikzset{>={Latex[width=3mm,length=2mm]}}
\tikzstyle{arr}=[->, thick, black]
\usetikzlibrary{decorations.text}
\usetikzlibrary{shadings,shapes.symbols}
\tikzset{
    double color fill/.code 2 args={
        \pgfdeclareverticalshading[%
            tikz@axis@top,tikz@axis@middle,tikz@axis@bottom%
        ]{diagonalfill}{100bp}{%
            color(0bp)=(tikz@axis@bottom);
            color(50bp)=(tikz@axis@bottom);
            color(50bp)=(tikz@axis@middle);
            color(50bp)=(tikz@axis@top);
            color(100bp)=(tikz@axis@top)
        }
        \tikzset{shade, left color=#1, right color=#2, shading=diagonalfill}
    }
}
\usepackage{subcaption}
\usepackage{amsfonts}
\usepackage{amssymb}
\usepackage{bbm}

\newcommand*{\qedb}{\null\nobreak\hfill\ensuremath{\square}}%

\newcommand{\ov}{\overline }

\def\A{\mathbf A}
\newcommand{\aaa}{\boldsymbol \alpha}
\newcommand{\bo}{\boldsymbol }
\newcommand{\nnu}{\boldsymbol\nu}

\newcommand{\yy}{\boldsymbol y}
\newcommand{\one}{\mathbf 1}
\newcommand{\zero}{\mathbf 0}
\newcommand{\B}{\mathbf B}
\newcommand{\I}{\mathbf I}
\newcommand{\E}{\mathbf E}
\newcommand{\F}{\mathbf F}

\newcommand{\G}{\mathbf G}
\newcommand{\X}{\mathbf X}
\newcommand{\Y}{\mathbf Y}
\def\P{\mathbf P}

\newcommand{\mt}{\mathcal T}
\newcommand{\MP}{\mathbb P}
\newcommand{\btheta}{\bo\Theta}

\newcommand{\bcolon}{\boldsymbol{:}}
\newcommand{\diag}{\text{\normalfont{diag}}}

\newcommand{\rank}{\text{\normalfont{rank}}}
\newcommand{\krank}{\text{\normalfont{krank}}}
\newcommand{\unfold}{\text{\normalfont{unfold}}}
\newcommand{\con}{\text{\normalfont{co}}}
\newcommand{\pa}{\text{\normalfont{pa}}}

\allowdisplaybreaks

\begin{document}

\title{Unfolding Tensors to Identify the Graph in Discrete Latent Bipartite Graphical Models}
\date{}
\author{Yuqi Gu\thanks{Email: \texttt{yuqi.gu@columbia.edu}}\\
Department of Statistics, Columbia University}

\maketitle

\begin{abstract}
We use a tensor unfolding technique to prove a new identifiability result for discrete bipartite graphical models, which have a bipartite graph between an observed and a latent layer. This model family includes popular models such as Noisy-Or Bayesian networks for medical diagnosis and Restricted Boltzmann Machines in machine learning. These models are also building blocks for deep generative models. Our result on identifying the graph structure enjoys the following nice properties. First, our identifiability proof is constructive, in which we innovatively unfold the population tensor under the model into matrices and inspect the rank properties of the resulting matrices to uncover the graph. This proof itself gives a population-level structure learning algorithm that outputs both the number of latent variables and the bipartite graph. Second, we allow various forms of nonlinear dependence among the variables, unlike many continuous latent variable graphical models that rely on linearity to show identifiability. Third, our identifiability condition is interpretable, only requiring each latent variable to connect to at least two ``pure'' observed variables in the bipartite graph. The new result not only brings novel advances in algebraic statistics but also has useful implications for these models' trustworthy applications in scientific disciplines and interpretable machine learning.
\end{abstract}

Tensor Unfolding.

\section{Introduction}

Probabilistic graphical models with latent variables are powerful tools to model complex joint distributions of random variables. They have been widely applied in statistics, machine learning, and various social and biomedical sciences \citep{lauritzen1996graphical, wainwright2008graphical, koller2009probabilistic}.
Incorporating unobserved {latent variables} 
in a graphical model often leads to very flexible and expressive generative models for high-dimensional heterogeneous data.
However, such flexibility comes at a cost of increasing model complexity, rendering the fundamental identifiability issue challenging to address.
Identifiability is a necessary prerequisite for valid statistical estimation and inference, and also an important step towards a model's interpretable and reproducible applications in scientific disciplines.
In this article, we use a {tensor unfolding} (flattening) technique to prove identifiability of the graph structure in discrete latent bipartite graphical models, which have a bipartite graph between one observed and one latent layer.

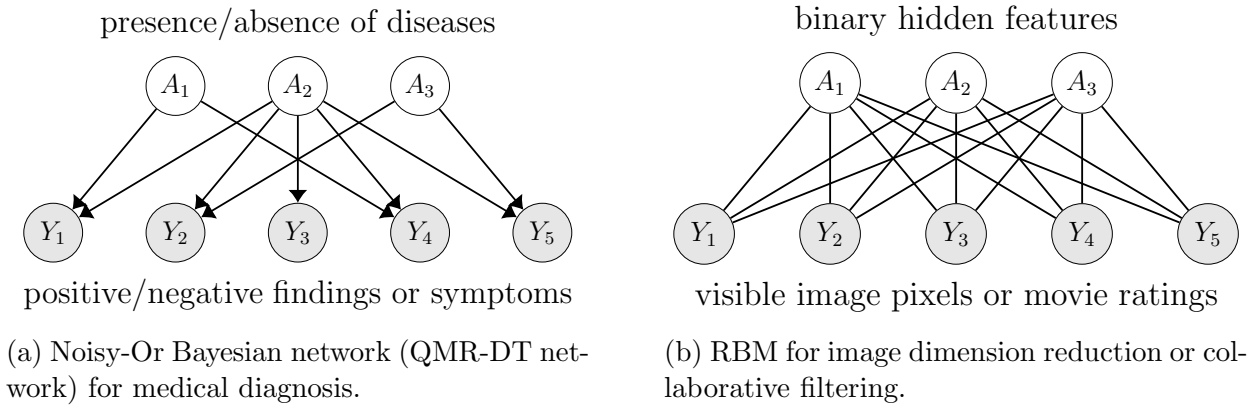
\begin{figure}[h!]\centering
\begin{subfigure}[b]{0.47\textwidth}
\resizebox{\textwidth}{!}{
\begin{tikzpicture}[scale=1.9]
    \node (v1)[neuron] at (0, 0) {$Y_{1}$};
    \node (v2)[neuron] at (1, 0) {$Y_{2}$};
    \node (v3)[neuron] at (2, 0) {$Y_{3}$};
    \node (v4)[neuron] at (3, 0) {$Y_{4}$};
    \node (v5)[neuron] at (4, 0) {$Y_{5}$};
       
    \node (h1)[hidden] at (1.0, 1.2) {$A_{1}$};
    \node (h2)[hidden] at (2.0, 1.2) {$A_{2}$};
    \node (h3)[hidden] at (3.0, 1.2) {$A_{3}$};
    
    \draw[->,thick,black] (h1) -- (v1) node [midway,above=-0.12cm,sloped] {\textcolor{black}{}}; 
    \draw[->,thick,black] (h1) -- (v4) node [midway,above=-0.12cm,sloped] {\textcolor{black}{}}; 

    \draw[->,thick,black] (h2) -- (v1) node [midway,above=-0.12cm,sloped] {\textcolor{black}{}}; 
    \draw[->,thick,black] (h2) -- (v2) node [midway,above=-0.12cm,sloped] {\textcolor{black}{}}; 
    \draw[->,thick,black] (h2) -- (v3) node [midway,above=-0.12cm,sloped] {\textcolor{black}{}}; 
    \draw[->,thick,black] (h2) -- (v4) node [midway,above=-0.12cm,sloped] {\textcolor{black}{}}; 
    \draw[->,thick,black] (h2) -- (v5) node [midway,above=-0.12cm,sloped] {\textcolor{black}{}}; 

    \draw[->,thick,black] (h3) -- (v2) node [midway,above=-0.12cm,sloped] {\textcolor{black}{}}; 
    \draw[->,thick,black] (h3) -- (v5) node [midway,above=-0.12cm,sloped] {\textcolor{black}{}};

     \node[anchor=south] (h2) at (2, 1.5) {\large \text{presence/absence of diseases}};

     \node[anchor=north] (h2) at (2, -0.3) {\large \text{positive/negative findings or symptoms}};
\end{tikzpicture}
}
\caption{Noisy-Or Bayesian network (QMR-DT network) for medical diagnosis.}
\end{subfigure}
\hfill
\begin{subfigure}[b]{0.47\textwidth}
\resizebox{\textwidth}{!}{
\begin{tikzpicture}[scale=1.9]
    \node (v1)[neuron] at (0, 0) {$Y_{1}$};
    \node (v2)[neuron] at (1, 0) {$Y_{2}$};
    \node (v3)[neuron] at (2, 0) {$Y_{3}$};
    \node (v4)[neuron] at (3, 0) {$Y_{4}$};
    \node (v5)[neuron] at (4, 0) {$Y_{5}$};
       
    \node (h1)[hidden] at (1.0, 1.2) {$A_{1}$};
    \node (h2)[hidden] at (2.0, 1.2) {$A_{2}$};
    \node (h3)[hidden] at (3.0, 1.2) {$A_{3}$};
    
    \draw[-,thick,black] (h1) -- (v1) node [midway,above=-0.12cm,sloped] {\textcolor{black}{}}; 
    \draw[-,thick,black] (h1) -- (v2) node [midway,above=-0.12cm,sloped] {\textcolor{black}{}}; 
    \draw[-,thick,black] (h1) -- (v3) node [midway,above=-0.12cm,sloped] {\textcolor{black}{}}; 
    \draw[-,thick,black] (h1) -- (v4) node [midway,above=-0.12cm,sloped] {\textcolor{black}{}}; 
    \draw[-,thick,black] (h1) -- (v5) node [midway,above=-0.12cm,sloped] {\textcolor{black}{}};  

    \draw[-,thick,black] (h2) -- (v1) node [midway,above=-0.12cm,sloped] {\textcolor{black}{}}; 
    \draw[-,thick,black] (h2) -- (v2) node [midway,above=-0.12cm,sloped] {\textcolor{black}{}}; 
    \draw[-,thick,black] (h2) -- (v3) node [midway,above=-0.12cm,sloped] {\textcolor{black}{}}; 
    \draw[-,thick,black] (h2) -- (v4) node [midway,above=-0.12cm,sloped] {\textcolor{black}{}}; 
    \draw[-,thick,black] (h2) -- (v5) node [midway,above=-0.12cm,sloped] {\textcolor{black}{}}; 

    \draw[-,thick,black] (h3) -- (v1) node [midway,above=-0.12cm,sloped] {\textcolor{black}{}}; 
    \draw[-,thick,black] (h3) -- (v2) node [midway,above=-0.12cm,sloped] {\textcolor{black}{}}; 
    \draw[-,thick,black] (h3) -- (v3) node [midway,above=-0.12cm,sloped] {\textcolor{black}{}}; 
    \draw[-,thick,black] (h3) -- (v4) node [midway,above=-0.12cm,sloped] {\textcolor{black}{}}; 
    \draw[-,thick,black] (h3) -- (v5) node [midway,above=-0.12cm,sloped] {\textcolor{black}{}}; 

     \node[anchor=south] (h2) at (2, 1.5) {\large \text{binary hidden features}};

     \node[anchor=north] (h2) at (2, -0.3) {\large \text{visible image pixels or movie ratings}};
\end{tikzpicture}
}
\caption{RBM for image dimension reduction or collaborative filtering.}
\end{subfigure}

\caption{Graphical model representations. Grey nodes denote observed variables, and white nodes denote latent variables.
\textbf{(a)}: QMR-DT network for medical diagnosis illustrated in \cite{jaakkola1999variational}. \textbf{(b)}: Restricted Boltzmann Machine (RBM) in \cite{hinton2006reducing} for dimension reduction of image data, and in \cite{salakhutdinov2007rbm} for collaborative filtering.}
\label{fig-exp1}
\end{figure}

Our considered model family includes two popular types of bipartite graphical models in statistics and machine learning: \emph{directed} graphical models such as Noisy-Or Bayesian networks \citep{pearl1988prob, shwe1991probabilistic}, and \emph{undirected} graphical models such as restricted Boltzmann machines \citep[RBMs;][]{hinton2006reducing} in machine learning.
Figure \ref{fig-exp1} gives the graphical model representations of a Noisy-Or network and a RBM, in which grey nodes denote observed variables and white nodes denote latent ones.
Noisy-Or Bayesian networks 
were used for medical diagnosis in the ``Quick Medical Reference'' (QMR-DT) network in \cite{shwe1991probabilistic}. 
In the medical diagnosis context, observed variables represent positive or negative symptoms exhibited by a patient, and latent variables represent this patient's hidden states of presence or absence of multiple diseases.
On the other hand,
RBMs were popularized by the influential \emph{Science} paper \cite{hinton2006reducing}, which proposed a scalable heuristic estimation algorithm for RBMs %called ``contrastive divergence'' 
and applied it to image data for dimension reduction.
Traditional RBMs have binary latent and observed variables.
\cite{salakhutdinov2007rbm} extended RBMs to tabular and count responses for collaborative filtering and recommendation systems.
In an RBM, observed variables represent the pixel values in an image or ratings of various movies given by a user, and latent variables represent the binary hidden features of the image or user.
RBMs are also applied in scientific disciplines including fMRI data analysis in neuroimaging \citep{hjelm2014rbmneural} and many-body quantum physics \citep{melko2019restricted}.

Directed and undirected graphical models are different machineries for encoding dependence among random variables.
In Section \ref{sec-id}, we will show that a universal sparsity condition on the bipartite graph is key to identifying the graph structure, regardless of whether the graph is directed or undirected. In short, our condition requires that each latent variable is connected to at least two ``pure'' observed variables, each of which is not connected to any other latent variable.
Hereafter, we will also call the bipartite graph the \emph{loading graph}, because it describes how the observed variables load on the latent ones.
The loading graph is a key quantity in multidimensional latent variable models.
Identifiability of the loading graph is crucial for unambiguously interpreting the latent variables and for ensuring the validity of statistical analyses in the considered models.

\subsection{Literature and Related Work}\label{subsec-lite}
The study of identifiability of the loading graph in \emph{continuous} latent variable models dates back to decades ago.
%Connections to previous identifiability proofs and results. 
The \emph{tetrad test} on the covariance matrix entries is used to identify the loading graph in Gaussian latent factor models \citep{glymour1987discovering, scheines1998tetrad}. 
See also \cite{silva2006latent} and \cite{bing2020overlap} for modern developments in related models. However, these covariance matrix-based identification approaches do not apply to the discrete variable models. The reason is that the joint distribution of multiple discrete variables takes the form of a higher-order tensor, which cannot be fully characterized by covariance matrices.  
For example, the joint distribution of three binary variables takes the form of a $2\times 2\times 2$ three-way tensor with entries given by their joint probability mass function.
Moreover, continuous latent variable models often impose stringent linearity assumptions to establish identifiability \citep{glymour1987discovering, silva2006latent, bing2020overlap}. 
It is largely unknown whether identifiability can be established when the model contains nonlinear transformations or interaction effects of variables, although incorporating these components usually greatly enhances a model's flexibility.

To identify discrete latent graphical models, the tensor structure must be taken into account. One powerful approach to proving identifiability in this regard is the Kruskal's theorem \citep{kruskal1977three}, which guarantees the uniqueness of three-way tensor decompositions under certain rank conditions.
\cite{allman2009} proposed a general framework to adapt Kruskal's theorem to prove identifiability in various discrete latent structure models.
Kruskal's theorem does not give a constructive approach to directly recovering any parameter from the population distribution.
Additionally, Kruskal's Theorem focuses on the uniqueness of the full tensor decomposition.
So, when it comes to identifying the arguably most interesting discrete graph structure, Kruskal's theorem may impose stronger than necessary conditions, as will be shown in Section \ref{sec-continuous}.

One interesting idea that partly motivated our tensor unfolding identification approach is the {quartet test} developed for discrete latent tree models (LTMs).
A LTM is a graphical model in which the graph is a tree and only the leaf variables are observed.
LTMs are popular in phylogenetics, statistics, and machine learning \citep{zwiernik2012tree, zwiernik2018chapter, shiers2016ltm, leung2018algebraic}.
Researchers have discovered that in discrete LTMs, the tensors corresponding to joint distributions of every four observed variables (i.e., quartets) can be unfolded to identify the latent tree structure \citep{eriksson2005tree, allman2006treetopology, ishteva2013unfolding, jaffe2021spectral}.
An important characteristic of LTMs is that they do not contain any cycle, implying that there is a unique path connecting any two variables. This fact turns out to be crucial for guaranteeing that all information about the graph is encoded in those fourth-order tensors for the quartets.
In contrast, the identifiability issue becomes more intriguing for complex graphical models beyond trees, such as our models. For example, there can be cycles in the bipartite graph: in Figure \ref{fig-exp1}(a),
the skeleton of the bipartite graph has a cycle $A_1 - Y_1 - A_2 - Y_4 - A_1$.

% \item
\subsection{Our Contributions}

% summarize our contributions
We propose a new identifiability condition for the bipartite graph in latent bipartite graphical models. Our proof is constructive and innovatively exploits tensor unfolding. Here is a high-level summary of the proof strategy.
Since all variables in the model are categorical, the joint distribution of the $J$ observed variables, upon marginalizing out the $K$ latent variables, is characterized by a $J$-way probability tensor. Entries of this tensor take complicated forms due to the marginalization -- they are sums over exponentially many hidden configurations of the $K$ latents.
For this $J$-way tensor, we carefully unfold it into matrices in various ways and then inspect the ranks of these matrices. Interestingly, we find that the bipartite graph induces remarkable rank properties of the unfolded matrices. These rank properties can serve as certificates for revealing the bipartite graph, as long as the true graph satisfies a simple condition: each latent variable is connected to at least two ``pure'' observed variables (see Theorem \ref{thm-main} for details).
Our constructive proof itself is indeed a population-level structure learning algorithm, which outputs both the number of latent variables and the bipartite graph.
% with the only input being the $J$-way tensor.
Remarkably, our identifiability result is {agnostic} and adaptive to whether the graph is directed or undirected, whether the latents are marginally dependent or independent, and whether there are nonlinear transformations in the model.

We point out that the latent bipartite graphical models can be much more flexible and complex than LTMs, and the new identifiability result is a significant advancement from existing results.
% significantly different from existing results for LTMs.
\emph{First}, the bipartite graph may contain cycles as mentioned earlier, rendering the quartet-based identification approach for LTMs not applicable.
\emph{Second}, by allowing a node to have multiple parents (resp. neighbors), the family of conditional distributions in the model is significantly enriched beyond tree-based models.
To see this, take directed graphs as an example. In a directed LTM, each variable has at most one parent under the tree structure.
The joint distribution of the LTM factorizes into the product of conditional distributions of each variable given its \emph{only parent}. 
In contrast, in our latent bipartite graphical model, each observed variable can have \emph{multiple} latent parents; e.g., in Figure \ref{fig-exp1}(a), $Y_1$ has two parents $A_1$ and $A_2$. As a result, the conditional distribution of a variable, given its multiple parents, has various flexible possibilities. One possibility is to incorporate both the main effects and interaction effects of those parents by letting $\MP(Y_1=1\mid A_1,A_2) = \beta_0+\beta_1 A_1 + \beta_2 A_2 + \beta_{12} A_1A_2$. Another possibility is to allow nonlinear transformations by letting $\MP(Y_1=1\mid A_1,A_2) = g(\beta_1 A_1 + \beta_2 A_2)$ with $g(x) = 1/(1+\exp(-x))$ or some other nonlinear function $g(\cdot)$.
So, establishing the identifiability of the bipartite graph for models in this generality, %as will be done in Theorem \ref{thm-main},
 is technically nontrivial and novel.

The new identifiability result not only brings technical advances, but also has useful practical implications.
On the one hand, when applying the considered models to social and biomedical sciences such as educational cognitive diagnosis of latent skills and medical diagnosis of latent diseases, our result will give an intuitive condition on the loading graph to help pinpoint the meaning of each latent trait. Indeed, the substantive meaning of each latent trait can be derived by looking at those ``pure'' observed features that are only connected to it.
On the other hand, in machine learning, the considered discrete bipartite graphical models have a huge representation power, as they induce a {factorial mixture} model with exponentially many mixture components, each corresponding to a configuration of the latent vector $\A$. 
{Recently, there has been a rapidly growing interest in establishing the identifiability for expressive machine learning models and deep generative models with latent variables \citep{khemakhem2020variational, moran2021identifiable,kivva2022identifiability,wu2023connecting,von2024nonparametric}. This is because only identifiable models can be learned reproducibly across training instances and deliver reliable representations for downstream tasks. In our considered models, establishing the bipartite graph's identifiability will help make these models more explainable and facilitate {interpretable machine learning}.}

%\vspace{4mm}
The rest of this article is organized as follows. Section \ref{sec-model} provides the general model setup and gives concrete examples.
Section \ref{sec-unfold} introduces the unfoldings of the population tensor.
Section \ref{sec-id} presents our main identifiability result along with an illustrative toy example.
Section \ref{sec-disc} includes concluding remarks. All of the proofs are included in the Supplementary Material.

\vspace{-5mm}
\section{General Model Setup}
\label{sec-model}

\subsection{Notations and Background}
For a set $S$, denotes its cardinality by $|S|$.
For a positive integer $M$, denote $[M]=\{1,2,\ldots,M\}$.
For a $J$-dimensional vector $\bo x = (x_1,\ldots,x_J)$ and any subset $S\subseteq[J]$, denote by $\bo x_S = (x_j)_{j\in S}$ the subvector of $\bo x$. %collecting all entries indexed by the integers in $S$.
For two vectors $\bo a = (a_1,\ldots,a_L)$ and $\bo b = (b_1,\ldots,b_L)$ of the same length, denote their entrywise maximum by $\bo a\vee\bo b = (\max(a_1,b_1),~\ldots,~ \max(a_L,b_L))$.
Denote by $\otimes$ the Kronecker product of matrices and by $\odot$ the Khatri-Rao product of matrices, which is the column-wise Kronecker product \citep{kolda2009tensor}.
In particular, consider matrices $\mathbf A=(a_{i,j})\in\mathbb R^{m\times r}$, $\mathbf B=(b_{i,j})\in\mathbb R^{s\times t}$; 
and matrices $\mathbf C=(c_{i,j})=(\bo c_{\bcolon,1}\mid\cdots\mid\bo c_{\bcolon,k})\in\mathbb R^{n\times k}$,
$\mathbf D=(d_{i,j})=(\bo d_{\bcolon,1}\mid\cdots\mid\bo d_{\bcolon,k})\in\mathbb R^{\ell\times k}$, then $\mathbf A\otimes \mathbf B \in\mathbb R^{ms\times rt}$ and $\mathbf C\odot \mathbf D \in\mathbb R^{n \ell\times k}$:
\begin{align*}
	\mathbf A\otimes \mathbf B
	=
	\begin{pmatrix}
		a_{1,1}\mathbf B & \cdots & a_{1,r}\mathbf B\\
		\vdots & \vdots & \vdots \\
		a_{m,1}\mathbf B & \cdots & a_{m,r}\mathbf B
	\end{pmatrix},
	\qquad
	\mathbf C\odot \mathbf D
	=
	\begin{pmatrix}
		\bo c_{\bcolon,1}\otimes\bo d_{\bcolon,1}
		\mid \cdots \mid
		\bo c_{\bcolon,k}\otimes\bo d_{\bcolon,k}
	\end{pmatrix}.
\end{align*}

Assume there are $K$ latent variables $A_1,\ldots,A_K$ and $J$ observed variables $Y_1,\ldots,Y_J$,
%in the graphical model, 
and a bipartite graph between them. 
Suppose each observed $Y_j$ ranges in $V$ categories $\{0,1,\ldots,V-1\}$ and each latent $A_k$ ranges in $H$ categories $\{0,1,\ldots,H-1\}$, where $V$, $H\geq 2$.
We use a $J \times K$ adjacency matrix $\G$ to describe the bipartite graph: $\G=(g_{j,k})\in\{0,1\}^{J\times K}$, where $g_{j,k} = 1$ or $0$ indicates the presence or absence of an edge between $Y_j$ and $A_k$. 
We allow for either directed or undirected edges. If the edges are directed, then we consider the case where all of them are pointing from the latent to the observed variables (as shown in Figure \ref{fig-exp1}(a)), following the convention of data generative processes in most latent variable models \citep{bishop2006pattern}.
Using the graph theory terminology, in a directed graphical model, $g_{j,k}=1$ means $A_k$ is a \emph{parent} of $Y_j$, whereas in an undirected graphical model,  $g_{j,k}=1$ means $A_k$ is a \emph{neighbor} of $Y_j$.
Define 
\begin{equation}\notag
    \con(j) = \{k\in[K]:~ g_{j,k}=1\},
\end{equation}
which collects the latent variables \textbf{co}nnected to $Y_j$ in the bipartite graph.
For a set $S\subseteq[p]$, let $\con(S) = \bigcup_{j\in S}\, \con(j)$ be the union of latent variables connected to the observed variables in $S$.

In general directed graphical models (DGMs) and undirected graphical models (UGMs), the presence of edges can imply different conditional independence relationships among variables in the graph.
% To see the above, one only needs to follow the graphical model definitions.
According to the definitions of probabilistic graphical models \citep[see, e.g.][]{wainwright2008graphical} that do not necessarily contain latent variables, the joint distributions of all variables $(X_1,\ldots,X_p)=:\X$ under a DGM and an UGM can be written as follows, respectively:
\begin{align}\label{eq-dgm}
\MP^{\text{DGM}}(\X) 
&= \prod_{j\in[p]} \MP(X_j\mid X_{\pa(j)});\quad
%\label{eq-ugm}
\MP^{\text{UGM}}(\X) 
= \frac{1}{Z} \prod_{C\in \mathcal C} \Psi_C(\X_{C}).
\end{align}
In \eqref{eq-dgm}, $\pa(j)\subseteq [p]\setminus\{j\}$ is the index set of all parent variables of $X_j$ in the DGM, and $\mathcal C$ is the collection of all maximal cliques in the UGM, $\Psi_C(\X_{C})$ is a positive function defined on the clique of variables indexed by $C\in\mathcal C$, and $Z$ is a constant chosen to ensure that the distribution is correctly normalized.
In the considered bipartite graphical model, if the edges in the bipartite graph are directed pointing from the latent layer to the observed layer, then the latent variables are marginally independent with each other. This fact can be directly seen from \eqref{eq-dgm}, because $\MP(\A) = \prod_{k=1}^K\MP(A_k)$ holds when there are no edges among the $K$ latent variables.
On the other hand, if the edges in the bipartite graph are undirected, then the latent variables are typically marginally dependent.
Despite this difference, in these two model families, the observed variables $Y_1,\ldots,Y_J$ are conditionally independent given the latent variables $A_1,\ldots,A_K$.
This property is called \emph{local independence}, which is an extremely common assumption in various latent variable models \citep{lazarsfeld1968latent}. 
Local independence implies that the latent variables account for how the observed responses are related to one another.

\subsection{Model Assumptions}
Since we consider categorical random variables in this work, 
all conditional distributions can be described by \emph{conditional probability tables} (CPTs). 
For notational simplicity, we use $\MP(Y_j\mid \A_{\con(j)})$ to denote the matrix or CPT that specifies the conditional distribution of $Y_j$ given its latent parents (resp. neighbors).
This matrix $\MP(Y_j\mid \A_{\con(j)})$ has size $V\times H^{|\con(j)|}$.
Its $H^{|\con(j)|}$ columns are indexed by all possible configurations of the latent vector $\A_{\con(j)}$ ranging in $\{0,1,\ldots,H-1\}^{|\con(j)|}$, and its $V$ rows are indexed by the $V$ categories of $Y_j$ ranging in $\{0,1,\ldots,V-1\}$.
%%%
We make the following two mild assumptions on $\MP(Y_j\mid \A_{\con(j)})$.

\begin{assumption}\label{assume-fullrank}
$V\geq H$. If $\con(j)$ is a singleton set, then the $V\times H$ conditional probability table $\MP(Y_j\mid A_{\con(j)})$ has full rank $H$.
\end{assumption}

\begin{assumption}
\label{assume-graph}
If $k\in\con(j)$, then there exists some vector $\bo a_{\con(j)\setminus\{k\}} \in \{0,1,\ldots,H-1\}^{|\con(j)\setminus\{k\}|}$ such that the following $H$ vectors are not identical:
    %are linearly independent:
    %\begin{align*}
        $\{\MP(Y_j\mid \A_{\con(j)\setminus\{k\}} = \bo a_{\con(j)\setminus\{k\}},~ A_k = h): ~~ h=0,1,\ldots,H-1\}$.
    %\end{align*}
    
\end{assumption}

Assumption \ref{assume-fullrank} requires that for those single-parent (or single-neighbor) variables $Y_j$, the $V\times H$ conditional probability table $\MP(Y_j\mid A_{\con(j)})$ has full column rank $H$. This assumption is common in the latent tree model literature, where each variable has exactly one parent. 
For example, 
when studying the latent tree model, \cite{ishteva2013unfolding}
assumed $V \geq H$ and all $V\times H$ conditional probability tables have full column rank.
Our considered models do not restrict each latent variable to have only one parent (or neighbor), and Assumption \ref{assume-fullrank} only requires those single-parent (or single-neighbor) observed variables to have full-column-rank CPTs.
For multi-parent (or multi-neighbor) observed variables in our models, the CPTs can have many more columns than rows, because $H^{|\con(j)|}$ can be quite large when $|\con(j)| \geq 2$. For these CPTs, we do not require them to have full rank and even allow them to have rank smaller than $H$.
% }\end{remark}

Assumption \ref{assume-graph} mainly concerns those multi-parent (or multi-neighbor) $Y_j$, because if $\con(j)$ is a singleton set, then Assumption \ref{assume-fullrank} implies Assumption \ref{assume-graph}. To see this, consider the case of $\con(j)=\{k\}$ for some $k\in[K]$. In this case, the $H$ vectors in Assumption \ref{assume-graph} are precisely the $H$ columns of the $V\times H$ matrix $\MP(Y_j\mid A_{\con(j)}) = \MP(Y_j\mid A_{k})$, so the full-rankness of this matrix required in Assumption \ref{assume-fullrank} implies that these $H$ column vectors are not identical, satisfying Assumption \ref{assume-graph}.

Assumption \ref{assume-graph} can be interpreted as imposing a very mild requirement that the model should somehow respect the graph structure.
In words, Assumption \ref{assume-graph} states that if $A_k$ is connected to $Y_j$ in the bipartite graph, then there must exist two latent vector configurations differing only in $A_k$ that will give rise to different conditional distributions of $Y_j\mid \A_{\con(j)}$.
Specifically, when each latent variable $A_k$ has exactly two categories with $H=2$, Assumption \ref{assume-graph} reduces to the following statement:
If $k\in\con(j)$, then there exists some vector $\bo a_{\con(j)\setminus\{k\}} \in \{0,1\}^{|\con(j)\setminus\{k\}|}$ such that
%\begin{align}\label{eq-as2-h2}
    $\MP(Y_j\mid \A_{\con(j)\setminus\{k\}} = \bo a_{\con(j)\setminus\{k\}},~ A_k = 0) \neq
    \MP(Y_j\mid \A_{\con(j)\setminus\{k\}} = \bo a_{\con(j)\setminus\{k\}},~ A_k = 1)$.
%\end{align}
Generally, Assumption \ref{assume-graph} is satisfied by a wide variety of linear and nonlinear, main-effect and interaction-effect models, including the examples in the next Section \ref{subsec-model}.

\subsection{Model Examples}\label{subsec-model}
We next present examples in two general classes: directed bipartite graphical models, and undirected ones. All these models can be shown to satisfy Assumptions \ref{assume-fullrank} and \ref{assume-graph} for generic parameters; here, following \cite{allman2009}, we say a statement holds ``for generic parameters'' or ``generically'', if the set of parameters for which the statement fails all reside in a Lebesgue measure zero subset of the parameter space.
Usually, these measure zero sets are simultaneous zero-sets of a finite collection of polynomials of model parameters, which are called \emph{algebraic varieties} \citep{allman2009}.
With this concept, Assumption \ref{assume-fullrank} is easy to satisfy because for those $j$ with $|\con(j)|=1$, the $V\times H$ matrix $\MP(Y_j\mid A_{\con(j)})$ has full rank for generic parameters when $V\geq H$.
As for Assumption \ref{assume-graph}, we will verify it when introducing the concrete examples.

\vspace{2mm}
\noindent
\textbf{First Family of Examples: Directed Bipartite Graphical Models.}
In these models, all the edges in the bipartite graph point from the latent layer to the observed layer. 
Following the general definition in \eqref{eq-dgm},
the marginal distribution of $\Y$ can be written as
\begin{align}\label{eq-gdgm}
    \MP(\Y = \yy) = \sum_{\bo a\in\{0,\ldots,H-1\}^K} \prod_{k=1}^K \MP(A_k = a_k) \prod_{j=1}^J \MP(Y_j = y_j\mid \A_{\con(j)} = \bo a_{\con(j)}),
\end{align}
for all $\yy \in\{0,\ldots,V-1\}^J$.
It remains to specify the conditional distributions $\MP(Y_j = y_j\mid \A_{\con(j)} = \bo a_{\con(j)})$.
For example, the popular Noisy-OR Bayesian Network \citep{shwe1991probabilistic} with binary $Y_j\in\{0,1\}$ and $A_k\in\{0,1\}$ has the following conditional probability:
\begin{align}\notag
%\label{eq-noisyor1}
    \MP^{\text{\normalfont{Noisy-OR}}}(Y_j = 0\mid \A_{\con(j)} = \bo a_{\con(j)}) 
    = \prod_{k:\; k\in\con(j)} \exp\left( -w_{j,k} a_k\right)
    = \exp\left( -\sum_{k:\; k\in\con(j)} w_{j,k} a_k \right).
    %\prod_{k:\; k\in\con(j)}
\end{align}
One canonical application of the above Noisy-OR Bayesian networks is to model how a patient's observable symptoms are related to the underlying disease statuses in medical diagnosis \citep{shwe1991probabilistic, jaakkola1999variational}. In this context, $V=H=2$ holds as $Y_j=1$ or 0 encodes the presence or absence of the $j$th symptom, and $A_k=1$ or 0 encodes the presence or absence of the $k$th latent disease. 
One can also define variations of the Noisy-OR Bayesian Network by using other link functions $f$ (such as the the identity link and the probit link):
\begin{align}\label{eq-noisyor2}
    \MP^{\text{\normalfont{Main-Effect}}}(Y_j = 0\mid \A_{\con(j)} = \bo a_{\con(j)}) 
    &=f\left(\sum_{k:\; k\in\con(j)} w_{j,k} a_k\right),
\end{align}
For the Noisy-OR Bayesian network and its variants, the number of parameters associated with the bipartite graph $\{w_{j,k}: \; k\in\con(j)\}$ equals the number of directed edges. All these networks model the main effects (or additive effects) of the  latent parents on the observed child variable.

To extend beyond the above main-effect models, we can define more complicated nonlinear models for $\Y\mid\A$ that include certain interaction effects of the latent parents.
Let $f(x)$ denote a monotone link function, which could be the identity link, logistic link, or probit link. Let $Y_j$ depend on all the possible main effects and interaction effects of its parent latent variables:
\begin{align}\label{eq-alleff}
    \MP^{\text{\normalfont{All-Effect}}}(Y_j = 1\mid \A_{\con(j)} = \bo a_{\con(j)})  = 
    f\left( \sum_{S\subseteq \con(j)} \beta_{j,S} \prod_{k\in S} a_k \right),
\end{align}
which has $2^{|\con(j)|}$ continuous parameters for each $j$. 
For example, if $\con(j) = \{A_1, A_3\}$, then $f\left( \sum_{S\subseteq \con(j)} \beta_{j,S} \prod_{k\in S} a_k \right) = f(\beta_{j,\varnothing} + \beta_{j,1}A_1 + \beta_{j,3}A_3 + \beta_{j,13}A_1 A_3)$ has four continuous parameters $\beta_{j,\varnothing}$, $\beta_{j,1}$, $\beta_{j,3}$, and $\beta_{j,13}$.
In the field of educational assessment in psychometrics, variants of the above model \eqref{eq-alleff} are very popular thanks to their flexibility \citep[e.g.][]{henson2009defining, von2008general, de2011generalized}. In those so-called \emph{cognitive diagnostic models}, one often has $V=H=2$: the latent variables encode presence or absence of a student's $K$ latent skills, the observed variables are the student's correct or wrong responses to $J$ test questions, and the bipartige graph $\G$ encodes which skills each test question is designed to measure.

As yet another example of an intermediate between the main-effect and all-effect models, one can let  $Y_j$ depend on both the main effects and the highest-order interaction effect of the latents:
\begin{align}\label{eq-mainhy}
    \MP^{\text{\normalfont{Main-Inter.}}}(Y_j = 1\mid \A_{\con(j)} = \bo a_{\con(j)})  = 
    f\left( \beta_{j,0} + \sum_{k\in\con(j)} \beta_{j,k} a_k +  \beta_{j,\text{all}} \prod_{k\in \con(j)} a_k \right),
\end{align}
which has $2+|\con(j)|$ parameters associated with the bipartite graph.
The term $\beta_{j,\text{all}} \prod_{k\in \con(j)} a_k$ captures the conjunctive relationship among the latent variables, meaning that mastering all of the required skills (i.e., parent latent skills) of a test question makes a difference on one's correct response probability. This is another very popular assumption in the educational cognitive diagnosis literature \citep{junker2001dina, gu2023joint}.
Now we show that Assumption \ref{assume-graph} is easy to satisfy for all of the above models. 
The conditional distribution of $Y_j\mid\A_{\con(j)}$ in \eqref{eq-gdgm}--\eqref{eq-mainhy} all involve the linear combinations of the latent parents/neighbors $\sum_{k\in\con(j)} \beta_{j,k}A_k$, so 
$\MP(Y_j\mid \A_{\con(j)}) = g\Big(\sum_{k\in\con(j)} \beta_{j,k}A_k + h( \A_{\con(j)})\Big)$,
 where $g(\cdot)$ and $h(\cdot)$ are potentially nonlinear functions.
Since $H=2$, if $k\in\con(j)$, then Assumption \ref{assume-graph} is satisfied as long as $\beta_{j,k}\neq 0$.

%%%%%%%%%%%%%%%%%%%%%%%%%%%%
\vspace{4mm}
\noindent
\textbf{Second Family of Examples: Undirected Bipartite Graphical Models.}
Restricted Boltzmann Machines are popular machine learning models with an undirected bipartite graph between an observed and a latent layer. 
RBMs are widely used in computer vision for dimension reduction of image data, such as facial images and images of handwritten digits \citep{hinton2006reducing}. They also serve as building blocks in deep generative models like deep belief networks \citep{goodfellow2016deep}.
In the most widely used binary RBM \citep{hinton2006reducing},
\begin{align}\notag
\MP^{\text{RBM}}(\Y = \yy) 
    &= \frac{1}{C}
     \sum_{\bo a\in\{0,1\}^K} \exp\left( \sum_{(j,k):\; g_{j,k}=1} \beta_{j,k} y_j a_k + \sum_{j=1}^J b_j y_j + \sum_{k=1}^K c_k a_k \right),~ \yy\in\{0,1\}^J.
\end{align}
In a general RBM \citep{salakhutdinov2007rbm} for categorical responses $\Y \in\{0,\ldots,V-1\}^J$ with categorical latent variables $\A\in \{0,\ldots,H-1\}^K$, we have
\begin{align}\label{eq-grbm}
& \MP^{\text{GenRBM}}(\Y = \yy) = 
\frac{1}{C}\sum_{\bo a\in\{0,\ldots,H-1\}^K}\exp\left( \sum_{(j,k):\; g_{j,k}=1} E_{jk}(y_j, a_k) + \sum_{j=1}^J E_j(y_j) + \sum_{k=1}^K E_k(a_k) \right),
\end{align}
where $C$ is the normalizing constant that ensures $\MP(\Y = \yy)$ is a valid probability mass function for $\Y$.
Specifically, since $y_j$ and $a_k$ are discrete, the terms $E_{jk}(y_j, a_k)$, $E_j(y_j)$, and $E_k(a_k)$ take the forms of 
$E_{jk}(y_j, a_k) = \sum_{v=0}^{V-1} \sum_{h=0}^{H-1} \beta_{j,k,v,h} \mathbbm{1}(y_j=v) \mathbbm{1}(a_k=h)$,
$E_j(y_j) = \sum_{v=0}^{V-1} b_{j,v} \mathbbm{1}(y_j=v)$,
$E_k(a_k) = \sum_{h=0}^{H-1} c_{k,h} \mathbbm{1}(a_k=h)$.
%\end{align*}
The latent variables in an RBM are generally \emph{not} marginally independent.
Now we verify that Assumption \ref{assume-graph} is easy to satisfy for \eqref{eq-grbm}. Fix any $j\in[J]$ and $k\in\con(j)$, the entries in the $h$th vector in the $H$ vectors in Assumption \ref{assume-graph} are
\begin{align*}
    &~\MP(Y_j=v\mid \A_{\con(j)\setminus\{k\}} = \bo a_{\con(j)\setminus\{k\}},~ A_k = h) \\
    \propto &~
    \exp\Big(b_{j,v} + \beta_{j,k,v,h} + c_{k,h} + \sum_{k'\in\con(j)\setminus\{k\}} \Big(E_{jk'}(v,a_{k'}) + E_{k'}(a_{k'})\Big) \Big).
\end{align*}
So, as long as the $H$ numbers $\{\beta_{j,k,v,h}:~ h=0,1,\ldots,H-1\}$ are not all equal for some $v$, 
the $H$ vectors in Assumption \ref{assume-graph} will not be identical generically and Assumption \ref{assume-graph} is satisfied.

\section{Unfoldings of the Population Tensor and Technical Preliminaries}
\label{sec-unfold}

Since all observed variables $Y_1,\ldots,Y_J$ are categorical, we can represent their joint distribution via a \emph{population probability tensor} $\mathcal T$ with $J$ ``modes'' and dimension $V$ along each mode. That is, the tensor $\mathcal T = (T_{i_1,i_2,\ldots,i_J})$ has size $V \times V \times \cdots \times V$, with entries being the following
$$
T_{i_1,i_2,\ldots,i_J} = \MP(Y_1=i_1,Y_2=i_2,\ldots,Y_J = i_J),\quad \forall i_1,i_2,\ldots,i_J \in[V].
$$
The tensor $\mathcal T$ can also be viewed as a $J$-way contingency table with $V^J$ cells, which is the conventional terminology for describing the joint distribution of multiple categorical random variables \citep{fienberg2009}.
For directed and undirected graphical models introduced in Section \ref{sec-model}, the marginal distribution of $\Y$ is given by \eqref{eq-gdgm} and \eqref{eq-grbm}, respectively.

\emph{Unfolding} (or flattening) a tensor refers to reshaping it into a matrix or a lower-order tensor by grouping certain modes together. 
For example, we can unfold a matrix $\mathbf C = (\bo c_1\mid \ldots\mid\bo c_p) \in\mathbb R^{p\times m}$ into a vector $(\bo c_1^\top, \ldots, \bo c_p^\top)^\top$ of dimension $pm$. 
As another example, for a 3rd-order tensor $\mathcal T = (T_{a,b,c}) \in \mathbb R^{d_1\times d_2 \times d_3}$, for any fixed $c \in[d_3]$, denote by $\mathcal T_{:, :, c}$ the $d_1\times d_2$ matrix whose $(a,b)$th entry is equal to $T_{a,b,c}$.
In this case, grouping the second and the third modes together results in a $d_1 \times d_2 d_3$ matrix 
$
({\mathcal T_{:, :, 1}}, \ldots, 
{\mathcal T_{:, :, d_3}}).
$
In our proof of the main result, it turns out that it suffices to consider unfolding the $J$-way population tensor $\mathcal T$ into various matrices (rather than other lower-order tensors). Each unfolded matrix comes from grouping certain modes in $[J]=\{1,\ldots,J\}$ into a row group and grouping the remaining modes in $[J]$ into a column group.
For example, if we consider the aforementioned unfolding of $\mathcal T = (T_{a,b,c}) \in \mathbb R^{d_1\times d_2 \times d_3}$ into a $d_1 \times d_2 d_3$ matrix, then the row group of modes is $\{1\}$ and the column group of modes is $\{2,3\}$.

Generally, define the \emph{row group} to be an arbitrary subset $S$ of $[J]$ and the \emph{column group} to be the remaining subset $[J]\setminus S$.
We write the tensor-unfolded matrix with row group $S$ and column group  $[J]\setminus S$ as ``$\unfold(\mathcal T,\; S,\; [J]\setminus S)$'' and use the following shorter notation for convenience:
$$
\unfold(\mathcal T,\; S,\; [J]\setminus S)
=: [\mathcal T]_{S, :}.
$$ 
Here, the unfolded matrix $[\mathcal T]_{S, :}$ has size $V^{|S|} \times V^{J-|S|}$, where the rows are indexed by all possible configurations of $\Y_S$ ranging in $\{0,\ldots,V-1\}^{|S|}$ and columns by all possible configurations of $\Y_{[J]\setminus S}$ ranging in $\{0,\ldots,V-1\}^{J-|S|}$.

Given an arbitrary subset $S\subseteq[J]$, define the \emph{marginal probability tensor} for $\Y_S$, which is a $|S|$-way tensor. The entries of this marginal probability tensor specify the joint probability mass function of the random vector $\Y_S$. It is easy to see that this marginal tensor can be obtained from the aforementioned full tensor $\mathcal T$ by appropriately summing up certain entries in it; for example, consider $\mathcal T = (T_{a,b,c}) \in \mathbb R^{d_1\times d_2 \times d_3}$, then the marginal tensor for $\Y_{\{1,2\}}$ is a $d_1\times d_2$ matrix, whose $(a,b)$-th entry is equal to $\sum_{c=1}^{d_3} T_{a,b,c}$.
So, $\mathcal T$ and $S$ uniquely defines the marginal probability tensor for $\Y_S$ and we denote it by 
$\text{marginal}(\mathcal T, S).$
Then, given two non-overlapping subsets $S_1$ and $S_2$ of $[J]$, the notation $\text{marginal}(\mathcal T, S_1\cup S_2)$ denotes the marginal probability tensor with $|S_1\cup S_2|$ modes. The following unfolding will also be useful in our proof:
$$
\unfold(\text{marginal}(\mathcal T, S_1\cup S_2),\; S_1,\; S_2) =:
[\mathcal T]_{S_1, S_2},
$$ 
which is a $V^{|S_1|} \times V^{|S_2|}$ matrix. This matrix essentially characterizes the joint distribution of $\Y_{S_1}$ and $\Y_{S_2}$. When $S_2 = [J]\setminus S_1$, it holds that $[\mathcal T]_{S_1, S_2} = [\mathcal T]_{S_1, :}$.

We introduce a useful notation of the \emph{joint probability table} between two random vectors with categorical entries.
Consider $A_1,\ldots,A_K\in\{0,1,\ldots,H-1\}$. 
For two sets $S_1, S_2\subseteq[K]$, we can describe the joint distribution of $\A_{S_1}=(A_k)_{k\in S_1}$ and $\A_{S_2}=(A_k)_{k\in S_2}$ using %a joint probability table 
$$
\MP(\A_{S_1}, \A_{S_2}): \text{ a $H^{|S_1|} \times H^{|S_2|}$ matrix}.
$$ 
The rows of $\MP(\A_{S_1}, \A_{S_2})$ are indexed by the configurations of $\A_{S_1}\in\{0,1,\ldots,H-1\}^{|S_1|}$, and columns by those of $\A_{S_2}\in\{0,1,\ldots,H-1\}^{|S_2|}$.
Importantly, we emphasize that $S_1$ and $S_2$ \emph{do not need to be disjoint and can overlap}. When $S_1\cap S_2\neq\varnothing$, $\MP(\A_{S_1}, \A_{S_2})$ has some zero entries corresponding to impossible configurations of $\A_{S_1\cup S_2}$. An extreme example is if $S_1=S_2$, then $\MP(\A_{S_1}, \A_{S_1})$ is a $H^{|S_1|} \times H^{|S_1|}$ diagonal matrix, because $\MP(\A_{S_1}=\bo a, \A_{S_1}=\bo b) = 0$ for any $\bo a\neq \bo b$. In this case, the diagonal entries of $\MP(\A_{S_1}, \A_{S_1})$ are given by the probability mass function of $\A_{S_1}$.
%%%
Another example is when $S_1 \subseteq S_2$, then the matrix $\MP(\A_{S_1}, \A_{S_2})$ has orthogonal row vectors because $\MP(\A_{S_1}=\bo a, \A_{S_2}=\bo b) \neq 0$ only if $\bo a$ is a subvector of $\bo b$ indexed by integers in $S_1$.
This general matrix notation of $\MP(\A_{S_1}, \A_{S_2})$ for potentially overlapping sets $S_1$ and $S_2$ turns out to be very useful to facilitate our identifiability proofs.
With this notation, we can also write the unfolded tensor $[\mathcal T]_{S,:}$ introduced in the previous paragraph as $\MP(\Y_{S}, \Y_{[J]\setminus S})$.

\section{Identifiability of the Bipartite Graph via Tensor Unfolding}\label{sec-id}

\subsection{Main Results}
Following the traditional definition of identifiability in the statistical literature \citep{koopmans1950identification, goodman1974}, we define the identifiability concept of the bipartite graph $\G$.

\begin{definition}\label{ref-defid}
The bipartite graph in the latent bipartite graphical model is said to be identifiable, if the $J\times K$ matrix $\G$ can be uniquely determined up to a column permutation from the joint distribution of the observed variables $\Y=(Y_1, \ldots, Y_J)$.
\end{definition}

In the considered models, the bipartite graph matrix $\G$ can at best be identified up to a permutation of its $K$ columns. Such permutation is an inevitable but trivial identifiability issue in related latent variable models; for example, in mixture models with $C$ latent classes, the mixture proportion parameters can at best be identified up to a label swapping of the $C$ latent classes \citep[e.g.,][]{allman2009}.
Therefore, as described in Definition \ref{ref-defid}, we focus on identifying $\G$ up to a permutation of the $K$ latent variables (equivalently, permutation of its $K$ columns).

The population distribution tensor $\mathcal T$ is the only available information to uncover the graph $\G$. Therefore, in order to prove identifiability under the mild Assumptions \ref{assume-fullrank} and \ref{assume-graph} for very general and flexible models described in Section \ref{sec-model}, we need to carefully investigate the algebraic structures hidden in $\mathcal T$ that can potentially imply the graph structure $\G$.
This is a technically nontrivial task, because all of the $K$ categorical latent variables have been marginalized out to obtain $\mathcal T$.
Indeed, each entry of $\mathcal T$ is a complicated function of the model parameters, as it is a sum of $H^K$ terms corresponding to all possible configurations of the latent vector $\A\in\{0,1,\ldots,H-1\}^K$. Each of these $H^K$ terms further involves polynomials and even nonlinear function transformations of model parameters in the considered models.

We have the following main theorem on the identifiability of $\G$, which is obtained via innovatively performing and examining the tensor unfoldings introduced in Section \ref{sec-unfold}.

\begin{theorem}\label{thm-main}
% Assume $V\geq H$.
Suppose Assumptions \ref{assume-fullrank} and \ref{assume-graph} hold and the $J\times K$ bipartite graph matrix $\G$ 
%(adjacency matrix for the bipartite graph between the latent and the observed variables) 
takes the following form after some row permutation:
% \singlespacing
\begin{align*}
    \G = \begin{pmatrix}
    \I_K\\
    \I_K\\
    \G^\star
    \end{pmatrix},
\end{align*}
where $\G^\star$ is an arbitrary binary matrix.
Then the bipartite graph $\G$ is identifiable. 
Moreover, the number of latent variables $K$ and the form of $\G$ can both be uniquely recovered using a constructive approach via repeatedly unfolding the $J$-way population distribution tensor $\mathcal T$.
%that characterizes the marginal distribution of $\Y$.
\end{theorem}

% discussing sparsity
The key identifiability condition in Theorem \ref{thm-main} is a sparsity requirement on the bipartite graph: each latent variable should be connected to at least two ``pure'' observed variables, each of which is only connected to this latent variable.
In scientific applications of the considered models, this condition translates to the requirement that each latent trait is measured by at least two observed features or items which solely measure this latent trait. Such a requirement is interpretable and shall be easily acceptable by applied researchers and practitioners.
Indeed, similar conditions have already been used in the educational assessment and cognitive diagnosis literature, where test designers require that for each latent skill, the exam should consist of some ``pure'' questions that solely measure this latent skill but not other skills \citep{von2019handbook}.

Next, we provide more technical details and insights into the new identifiability result and its proof.
Specifically, we prove Theorem \ref{thm-main} by establishing two key propositions, the following Propositions \ref{prop1} and \ref{prop2}, which reveal the single-parent/neighbor structures and the multiple-parent/neighbor structures in the bipartite graph, respectively.

\begin{proposition}[Unfold $\mathcal T$ to Reveal All Single-parent/neighbor Structures]
\label{prop1}
Under the condition in Theorem \ref{thm-main}, for arbitrary $j_1\neq j_2\in[J]$, it holds that 
$$\rank([\mt]_{\{j_1,j_2\},\; :}) \leq H
\quad\text{if and only if}\quad
\con(j_1) = \con(j_2) = \{k\}\text{  for some } k\in[K].$$ 
\end{proposition}

An important and remarkable property of our tensor unfolding identification approach is that, the number of latent variables $K$ can also be constructively identified from the tensor $\mathcal T$. The reason behind this can be directly deduced from Proposition \ref{prop1}. The statement in Proposition \ref{prop1}, if proved to be true, implies that one can find $K^*$ subsets of $[J]$ denoted by $S_1,\ldots,S_{K^*}$, which exhibit the following property:
\begin{align*}
    &\forall k_1\neq k_2\in[K^*],\quad S_{k_1} \cap S_{k_2} = \varnothing;
    \quad\forall k\in[K^*],
    \quad \forall j_1\neq j_2\in S_k,
    \quad \rank([\mt]_{\{j_1,j_2\},\; :}) \leq H.
\end{align*}
In words, the above display means that there exist $K^*$ mutually exclusive subsets of the $J$ observed variables, such that if we unfold the tensor $\mathcal T$ with the row group consisting of any two variables from any one of these subsets, the rank of this resulting matrix is no greater than $H$.
It is not hard to see that the maximum number that $K^*$ can take is precisely $K$, the true number of latent variables. Indeed, for $k\in[K]$, the set $S_k$ is characterized by 
$S_k=\{j\in[J]:~ \con(j)=\{k\}\},$
so $S_k$ collects the indices of all those pure observed variables only connected to $A_k$.
This nice property of constructively identifying the number of latents is in sharp contrast to previous identifiability results for related models, which were proved by using Kruskal's Theorem \citep[e.g., in][]{allman2009, xu2017rlcm, culpepper2019ordinal, gu2023bp}. Adapting Kruskal's Theorem to show identifiability always requires the knowledge of the number of latent variables.

Proposition \ref{prop1} provides a ``certificate'' to uncover the number of latent variables and all those single-parent (resp. single-neighbor) structures in the graph. In other words, the number of columns in the matrix $\G$ as well as all its standard basis row vectors have been identified and recovered from $\mathcal T$. 
Now, it remains to recover the remaining row vectors of $\G$, or equivalently, to uncover the graph structures for those observed variables who have more than one parents (resp. neighbors). 
This is what the following Proposition \ref{prop2} establishes.

\begin{proposition}[Unfold $\mathcal T$ to Reveal All Multi-parent/neighbor Structures]
\label{prop2}
Under the condition in Theorem \ref{thm-main}, assume without loss of generality that the first $2K$ rows of $\G$ are identified to be $(\I_K;~ \I_K)^\top$.
The following holds for arbitrary $k\in[K]$ and $j\in \{2K+1,\ldots,J\}$ with $|\con(j)|\geq 2$:
\begin{align}\label{eq-setab}
\rank([\mathcal T]_{([K]\setminus\{k\}) \cup \{j\},~ \{K+1,\ldots,2K\}})> H^{K-1}
\quad\text{if and only if}\quad k\in \con(j).
\end{align}
\end{proposition}

Proposition \ref{prop2} establishes another certificate for checking whether an arbitrary latent $A_k$ is connected to an arbitrary observed $Y_j$ which has at least two latent parents/neighbors.
In terms of the $\G$ matrix, Proposition \ref{prop2} hence identifies all of the remaining row vectors in $\G$, which are indexed by those $j\in[J]$ for which $|\con(j)| \geq 2$. 
Now, it is clear that combining Proposition \ref{prop1} and Proposition \ref{prop2} delivers the desired constructive identification result in Theorem \ref{thm-main}.

Also, it is possible to avoid accessing the whole $J$-way tensor but rather just unfold certain lower-order marginal tensors to show identifiability.
The next corollary formalizes this statement.

\begin{corollary}\label{cor-2k}
Under the condition in Theorem \ref{thm-main},
$\G$ is identifiable from unfolding the collection of marginal tensors of the following form:
$\{\text{marginal}(\mathcal T, S),~ \forall S\subseteq[J]\text{ and } |S|=2K.\}$
\end{corollary}

\subsection{A Toy Example Illustrating the Proof Idea}
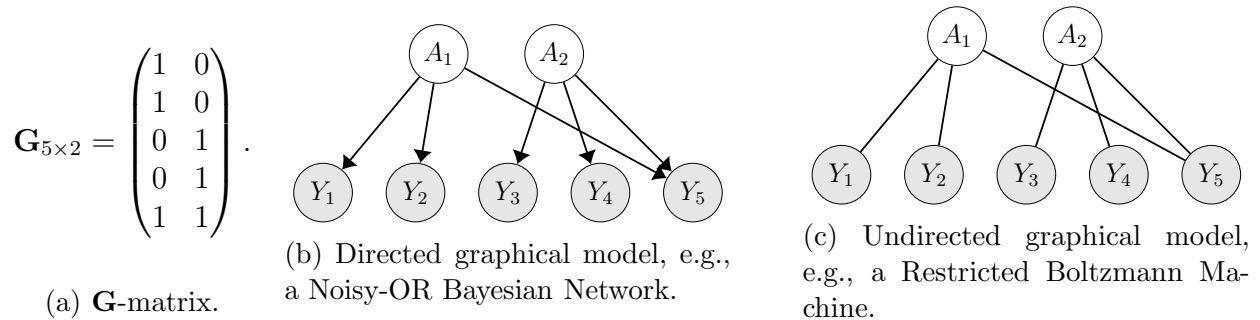
\begin{figure}[h!]\centering
\begin{subfigure}[c]{0.2\textwidth}
$$\G_{5\times 2} =
\begin{pmatrix}
    1 & 0 \\
    1 & 0 \\
    0 & 1 \\
    0 & 1 \\
    1 & 1
\end{pmatrix}.$$
\caption{$\G$-matrix.}
\end{subfigure}
~
\begin{subfigure}[c]{0.36\textwidth}
\resizebox{\textwidth}{!}{
\begin{tikzpicture}[scale=1.8]
    \node (v1)[neuron] at (0, 0) {$Y_{1}$};
    \node (v2)[neuron] at (0.8, 0) {$Y_{2}$};
    \node (v3)[neuron] at (1.6, 0) {$Y_{3}$};
    \node (v4)[neuron] at (2.4, 0) {$Y_{4}$};
    \node (v5)[neuron] at (3.2, 0) {$Y_{5}$};
       
    \node (h1)[hidden] at (1.0, 1.2) {$A_{1}$};
    \node (h2)[hidden] at (2.0, 1.2) {$A_{2}$};

    \draw[qedge] (h1) -- (v1) node [midway,above=-0.12cm,sloped] {\textcolor{black}{}}; 
    
    \draw[qedge] (h1) -- (v2) node [midway,above=-0.12cm,sloped] {};  
    
    \draw[qedge] (h2) -- (v3) node [midway,above=-0.12cm,sloped] {}; 
    
    \draw[qedge] (h2) -- (v4) node [midway,above=-0.12cm,sloped] {}; 
    
    \draw[qedge] (h1) -- (v5) node [midway,above=-0.12cm,sloped] {}; 
    
    \draw[qedge] (h2) -- (v5) node [midway,above=-0.12cm,sloped] {}; 
\end{tikzpicture}
}
\caption{Directed graphical model, e.g., a Noisy-OR Bayesian Network.}\end{subfigure}
\qquad
\begin{subfigure}[c]{0.36\textwidth}
\resizebox{\textwidth}{!}{
\begin{tikzpicture}[scale=1.8]
    \node (v1)[neuron] at (0, 0) {$Y_{1}$};
    \node (v2)[neuron] at (0.8, 0) {$Y_{2}$};
    \node (v3)[neuron] at (1.6, 0) {$Y_{3}$};
    \node (v4)[neuron] at (2.4, 0) {$Y_{4}$};
    \node (v5)[neuron] at (3.2, 0) {$Y_{5}$};
       
    \node (h1)[hidden] at (1.0, 1.2) {$A_{1}$};
    \node (h2)[hidden] at (2.0, 1.2) {$A_{2}$};
    
    \draw[-,thick,black] (h1) -- (v1) node [midway,above=-0.12cm,sloped] {\textcolor{black}{}}; 
    
    \draw[-,thick,black] (h1) -- (v2) node [midway,above=-0.12cm,sloped] {};  
    
    \draw[-,thick,black] (h2) -- (v3) node [midway,above=-0.12cm,sloped] {}; 
    
    \draw[-,thick,black] (h2) -- (v4) node [midway,above=-0.12cm,sloped] {}; 
    
    \draw[-,thick,black] (h1) -- (v5) node [midway,above=-0.12cm,sloped] {}; 
    
    \draw[-,thick,black] (h2) -- (v5) node [midway,above=-0.12cm,sloped] {}; 
\end{tikzpicture}
}
\caption{Undirected graphical model, e.g., a Restricted Boltzmann Machine.}
\end{subfigure}
\caption{Directed and undirected bipartite graphical models with a latent layer $\A=(A_1,A_2) \in [H]^2$ and an observed layer $\Y = (Y_1,\ldots,Y_5) \in [V]^5$.
$A_1$ and $A_2$ are \emph{marginally independent} in (a) and \emph{marginally dependent} in (b). Theorem \ref{thm-main} applies to both cases.}
\label{fig-K2}
\end{figure}

We present a toy example with $(J,K)=(5,2)$ to illustrate the tensor unfoldings under the identifiability condition in Theorem \ref{thm-main}.
The population tensor $\mathcal T$ has size $V \times V \times V \times V \times V$. Consider $\G_{5\times 2}$ specified in Figure \ref{fig-K2}(a), which satisfies the identifiability condition in Theorem \ref{thm-main} because its rows can be swapped so that it contains two $\I_2$.
Next, we will work through various unfoldings of $\mathcal T$ to reveal the remarkable rank properties the bipartite graph leaves on them.

Recall the joint probability table $\MP(\A_{S_1}, \A_{S_2})$ introduced in Section \ref{sec-unfold} and consider $\MP(A_1, \A_{1,2})$.
% representing the joint distribution of $A_1$ and $\A_{1,2}$.
Since $A_1$ is an entry of $\A_{1,2}$, the discussion in Section \ref{sec-unfold} implies that the $H \times H^2$ matrix $\MP(A_1, \A_{1,2})$ has orthogonal rows with mutually exclusive support.
For any $a_1\in\{0,\ldots,H-1\}$, $(b_1,b_2)\in \{0,\ldots,H-1\}^2$, 
\begin{align*}
    \MP(A_1, \A_{1,2})_{b,\; (a_1,a_2)} 
    &= \MP(A_1 = b,\; \A_{1,2} = (a_1,a_2))
    = 
    \begin{cases}
    0, & \text{if } b \neq a_1;\\
    \MP(\A_{1,2} = (a_1,a_2)), & \text{if } b = a_1.
    \end{cases}
\end{align*}
In the simpler special case with $H=2$, the $2\times 4$ matrix $\MP(A_1, \A_{1,2})$ can be written as:
\begin{align*}
    \MP(A_1, \A_{1,2}) = 
    {\small\begin{pmatrix}
        \MP(A_1=0,A_2=0) &~~ 0 &~~ \MP(A_1=0,A_2=1) &~~ 0 \\
        0 &~~ \MP(A_1=1,A_2=0) &~~  0 &~~\MP(A_1=1,A_2=1)
    \end{pmatrix}},
\end{align*}
where the zero entries indicate that the corresponding row configuration of $A_1$ and column configuration of $\A_{1,2}$ are not compatible. The first useful observation is that the rank of $\MP(A_1, \A_{1,2})$ equals $H$ generically. Specifically, the matrix $\MP(A_1, \A_{1,2})$ has full rank as long as $\MP(A_1 = 0) \in (0,1)$, because the row sums of this nonnegative matrix are $\MP(A_1 = 0)$ and $\MP(A_1 = 1)$, respectively, and full-rankness holds if both $\MP(A_1 = 0)$ and $\MP(A_1 = 1)$ are nonzero.

We unfold the tensor $\mathcal T$ by defining the row group to be $\{1,2\}$ and the column group to be $\{3,4,5\}$, so the resulting unfolding is the joint probability table $\mathbb P(\Y_{1,2},\; \Y_{3,4,5})$. 
To understand the rank property of $\mathbb P(\Y_{1,2},\; \Y_{3,4,5})$, first note that $\con(1)=\con(2)=\{1\}$ in the bipartite graph (either directed or undirected) implies the following for any fixed vectors $\bo y_{1,2}$ and $\bo y_{3,4,5}$:
\begin{align}\label{eq-a12b}
    &~\mathbb P(\Y_{1,2}=\bo y_{1,2},\; \Y_{3,4,5}=\bo y_{3,4,5})\\ \notag
    =&~
    \sum_{\bo a_{1,2}\in\{0,\ldots,H-1\}^2} \MP(\A_{1,2}=\bo a_{1,2}) 
    \prod_{j=1}^5
    \MP(Y_j=y_j\mid \A_{\con(j)}=\bo a_{\con(j)}) \\ \notag
    %=&~ \sum_{\bo a\in\{0,\ldots,H-1\}^2} \MP(\A=\bo a)\MP(Y_1=y_1\mid A_1=a_1)  \MP(Y_2=y_2\mid A_1=a_1) \MP(\Y_{3,4,5} = \bo y_{3,4,5}) \\
    =&~ \sum_{\bo a_{1,2}\in\{0,\ldots,H-1\}^2} \MP(\A_{1,2}=\bo a_{1,2})
    \MP(\Y_{1,2}=\bo y_{1,2}\mid A_1=a_1)
    \MP(\Y_{3,4,5} = \bo y_{3,4,5}\mid \A_{1,2}=\bo a_{1,2}) \\ \notag
    =&~ \sum_{\bo a_{1,2}\in\{0,\ldots,H-1\}^2} \sum_{b\in\{0,\ldots,H-1\}}
    \MP(\Y_{1,2}=\bo y_{1,2}\mid A_1=a_1)
    \MP(A_1=b, \A_{1,2}=\bo a_{1,2})\\ 
    \notag
    &\qquad\qquad\qquad\qquad\qquad\qquad \times
    \MP(\Y_{3,4,5} = \bo y_{3,4,5}\mid \A_{1,2}=\bo a_{1,2}),
\end{align}
where $b\in \{0,\ldots,H-1\}$ in the last equality above is deliberately introduced in order to further write $[\mt]_{\{1,2\},\; :}$ as a matrix factorization. 
Based on the above equation, if we let the vectors $\bo y_{1,2}$ and $\bo y_{3,4,5}$ range in $\{0,\ldots,V-1\}^2$ and $\{0,\ldots,V-1\}^3$, respectively, then we can obtain the following key factorization of the $V^2 \times V^3$ matrix $[\mt]_{\{1,2\},\; :}$
resulting from the tensor unfolding: 
\begin{align}\label{eq-toyrank}
    [\mt]_{\{1,2\},\; :} 
    = 
    \left( \mathbb P(\Y_{1,2},\; \Y_{3,4,5})
    \right)_{V^2 \times V^3}
    &\stackrel{(\star)}{=} 
    \underbrace{\mathbb P(\Y_{1,2}\mid A_1)}_{V^2 \times H} 
    \cdot \underbrace{\mathbb P(A_1, \A_{1,2})}_{H\times H^2}
    \cdot \underbrace{\mathbb P(\Y_{3,4,5}\mid \A_{1,2})^\top}_{H^2\times V^3},
    \\ \notag
    &
    \Longrightarrow\quad
    \rank([\mt]_{\{1,2\},\; :} ) \leq H.
\end{align}
The key equality marked with $(\star)$ above utilizes equation \eqref{eq-a12b} and the previously introduced notation of the joint probability table $\mathbb P(A_1, \A_{1,2})$.
According to the decomposition in \eqref{eq-toyrank}, the $4\times 8$ matrix $[\mt]_{\{1,2\},\; :}$ has rank at most $H$. 
Similarly, since $\con(3)=\con(4)=\{2\}$, we can use the same reasoning as above to obtain
\begin{align}\notag
    [\mt]_{\{3,4\},\; :} 
    &= 
    \underbrace{\mathbb P(\Y_{3,4}\mid A_2)}_{V^2 \times H} 
    \cdot \underbrace{\mathbb P(A_2, \A_{1,2})}_{H\times H^2}
    \cdot \underbrace{\mathbb P(\Y_{1,2,5}\mid \A_{1,2})^\top}_{H^2 \times V^3},\\ \label{eq-t34}
    &
    \Longrightarrow\quad \rank([\mt]_{\{3,4\},\; :} ) \leq H.
\end{align}
On the other hand, consider unfolding $\mt$ with the row group $\{1,3\}$ and column group $\{2,4,5\}$, and we can obtain the factorization using Lemma \ref{lem-redun-latent} in the Appendix:
\begin{align}\label{eq-t13exp}
    [\mt]_{\{1,3\},\; :} 
    &= \underbrace{\mathbb P(\Y_{1,3}\mid \A_{1,2})}_{\mathbf P_1:~V^2 \times H^2} 
    \cdot \underbrace{
    %\diag(\mathbb P(\A_{1,2}))
    \MP(\A_{1,2}, \A_{1,2})
    }_{\mathbf P_2:~H^2\times H^2}
    \cdot \underbrace{\mathbb P(\Y_{2,4,5}\mid \A_{1,2})^\top}_{\mathbf P_3:~H^2\times V^3}.
\end{align}
Next, we will show that the rank of the above matrix $[\mt]_{\{1,3\},\; :}$ is \emph{strictly greater than} $H$, distinguishing it from $[\mt]_{\{1,2\},\; :}$ and $[\mt]_{\{2,4\},\; :}$ examined earlier.
As long as each configuration of $\A_{1,2}$ has a nonzero probability with $\MP(\A_{1,2}=\bo a_{1,2}) > 0$ for any $\bo a_{1,2}\in\{0,\ldots,H-1\}^2$, we have $\rank(\mathbf P_2)=\rank(\MP(\A_{1,2}, \A_{1,2})) = H^2$.
% Since $\nb_{\G}(1)=1$ and $\nb_{\G}(3)=2$, we have 
Since $\con(1)=\con(2)=\{1\}$, $\con(3)=\con(4)=\{2\}$, and $\con(5)=\{1,2\}$, we apply Lemma \ref{lem-redun-latent} in the Appendix to write: 
\begin{align*}
\mathbf P_1=\underbrace{\mathbb P(\Y_{1,3}\mid \A_{1,2})}_{V^2 \times H^2} 
&= \underbrace{\MP(Y_1\mid A_1)}_{V\times H} \otimes \underbrace{\MP(Y_3\mid A_2)}_{V\times H},\\
\mathbf P_3=\underbrace{\mathbb P(\Y_{2,4,5}\mid \A_{1,2})}_{V^3 \times H^2} 
&=
\underbrace{\left(\MP(Y_2\mid A_1) \otimes \MP(Y_4\mid A_2)\right)}_{V^2\times H^2} \odot \underbrace{\MP(Y_5\mid \A_{1,2})}_{V\times H^2};
\end{align*}
In the above display, we slightly abuse notations when writing the Kronecker product and Khatri-Rao product, in that the equation holds up to potential permutations of the columns and rows in the product. Such permutations will not impact any of our proof arguments and final conclusion. This is because to establish identifiability, we only consider the matrix ranks of the Kronecker product and Khatri-Rao product, which are not affected by column and row permutations.
By Assumption \ref{assume-fullrank}, both $\MP(Y_1\mid A_1)$ and $\MP(Y_3\mid A_2)$ have full rank $H$, so their Kronecker product $\mathbb P(\Y_{1,3}\mid \A_{1,2})$ has full rank $H^2$. By a similar reasoning, $\MP(Y_2\mid A_1) \otimes \MP(Y_4\mid A_2)$ also has full column rank $H^2$. We further use Lemma \ref{lem-khatri}(b) in the Supplementary Material to obtain 
$$
\rank(\mathbf P_3) = \rank(\MP(\Y_{2,4,5}\mid \A_{1,2})) \geq \rank(\MP(Y_2\mid A_1) \otimes \MP(Y_4\mid A_2)) =
H^2.
$$
Going back to $[\mt]_{\{1,3\},\; :}=\mathbf P_1 \mathbf P_2 \mathbf P_3$ in \eqref{eq-t13exp}, since all three matrices have full rank, we have
\begin{align}\notag
&\rank([\mt]_{\{1,3\},\; :}) =H^2 > H. \\
\label{eq-t13exp}
&(\text{By symmetry, }
\rank([\mt]_{\{1,4\},\; :}), \rank([\mt]_{\{2,3\},\; :}), \rank([\mt]_{\{2,4\},\; :})
> H \text{ also hold.})
\end{align}
Last, consider $[\mt]_{\{1,5\},\; :}$ where $\con(1)=\{1\}$ and $\con(5)=\{1,2\}$:
\begin{align*}
    [\mt]_{\{1,5\},\; :} = 
    \underbrace{\MP(\Y_{1,5}\mid \A_{1,2})}_{\mathbf P_4:~V^2 \times H^2} \cdot \underbrace{\MP(\A_{1,2}, \A_{1,2})}_{\mathbf P_5:~H^2 \times H^2} \cdot \underbrace{\MP(\Y_{2,3,4}\mid \A_{1,2})^\top}_{\mathbf P_6:~H^2 \times V^3},
\end{align*}
where $\mathbf P_5$ is equal to $\mathbf P_2$ defined earlier, so $\rank(\mathbf P_5) = H^2$ for generic parameters. 
As for $\mathbf P_4$, we use Lemma \ref{lem-j1j2} in the Supplementary Material to obtain $\rank(\mathbf P_4)>H$.
As for $\mathbf P_6$, we can use a similar argument as showing $\rank(\mathbf P_3)=H^2$ to show $\rank(\mathbf P_6) = H^2$. 
Now we have shown that $[\mt]_{\{1,5\},\; :}$ is the product of three matrices $\mathbf P_4, \mathbf P_5, \mathbf P_6$, with two of them ($\mathbf P_5$ and $\mathbf P_6$) having full rank $H^2$ and the other one ($\mathbf P_4$) having rank greater than $H$. As a result, we obtain
\begin{align}\notag
&\rank([\mt]_{\{1,5\},\; :}) = \rank(\mathbf P_4) > H.\\
\label{eq-t15exp}
&(\text{By symmetry, }
\rank([\mt]_{\{2,5\},\; :}), \rank([\mt]_{\{3,5\},\; :}), \rank([\mt]_{\{4,5\},\; :}) > H \text{ also hold.})
\end{align}
Summarizing \eqref{eq-toyrank}, \eqref{eq-t34}, \eqref{eq-t13exp}, and \eqref{eq-t15exp}, when we examine the rank of each of the 10 unfoldings $[\mt]_{\{j_1,j_2\},\; :}$ for all pairs of $j_1, j_2\in\{1,2,3,4,5\}$, the only two unfoldings that have rank no greater than $H$ are $[\mt]_{\{1,2\},\; :}$ and $[\mt]_{\{3,4\},\; :}$, which exactly correspond to the pure connection structures $\con(1)=\con(2)=\{1\}$ and $\con(3)=\con(4)=\{2\}$.
This rank property can be used to reveal the single-parent/-neighbor structures by assigning $Y_1$ and $Y_2$ to be the pure connections of one latent variable and $Y_3$ and $Y_4$ to be the pure connections of another one. The number of latent variables is also correctly identified as two here.
This delivers a very important insight: unfolding the tensor into a matrix with the row group containing exactly two observed variables reveals the structure of all pure children/neighbors.
We have now illustrated the key idea behind the proof of Proposition \ref{prop1}, although proving the fully general case requires more technical effort.

Proving Proposition \ref{prop2} for revealing the multi-parent/neighbor structures is more difficult and hence not illustrated in great detail here, but the core idea in proving Proposition \ref{prop2} is still to examine the ranks of the unfoldings. The difficulties include: first, to come up with an appropriate row group $S_1$ and a corresponding column group $S_2$ of variables from $Y_1,\ldots, Y_J$ such that $[\mathcal T]_{S_1,S_2}$ contains as much as possible useful information about the graph; and second, to carefully investigate the rank of these unfolded matrices to establish the certificate for the graph.

\subsection{Identifiability of Other Parameters in the Model}
\label{sec-continuous}
Next, we present a side result on the identifiability of other continuous parameters in the considered model and explain how it is connected to the main result of identifying $\G$. 
We first introduce a unified notation for the continuous parameters. For any $\bo a\in\{0,\ldots,H-1\}^K$ and any $v\in\{0,1,\ldots,V-1\}$, define
    $\nu_{\bo a} = \mathbb P(\A = \bo a)$, and
    %\qquad
    $\theta^{(j)}_{v\mid \bo a_{\pa(j)}} =
    \mathbb P(Y_j=v\mid\A_{\pa(j)} = \bo a_{\pa(j)})$.
%\end{align*}
Let $\nnu$ collects all the $\nu_{\bo a}$-parameters and $\bo\Theta$ collects all the $\theta^{(j)}_{v\mid \bo a_{\pa(j)}}$-parameters. 
The definition of $\bo\Theta$ reflects the influence of the graph $\G$ on the conditional distribution $Y\mid\A$. The parameters $\nnu$ and $\bo\Theta$ fully characterize the joint distribution of $\A$ and $\Y$. 
Consider the following parameter space for $\nnu$ and $\bo\Theta$:
\begin{align*}
    \nnu\in & \Big\{\nu_{\aaa} > 0,~\sum_{\aaa\in\{0,\ldots,H-1\}^K} \nu_{\aaa} = 1\Big\};\quad
    \bo\Theta\in \Big\{\theta^{(j)}_{v\mid \bo a_{\pa(j)}} > 0,\quad \sum_{v=0}^{V-1}\theta^{(j)}_{v\mid \bo a_{\pa(j)}} = 1, \quad\forall \bo a_{\pa(j)}%\in \{0,\ldots,H-1\}^{|\pa(j)|}
    \Big\}.
\end{align*}
The following proposition describes how the graph $\G$ impacts the identifiability of $\nnu$ and $\bo\Theta$.
%$\nnu = (\nu_{\aaa}: \aaa\in\{)$

\begin{proposition}\label{prop-cont}
%[Theorem 3.3 in \cite{gu2024blessing}] 
The following conclusions hold.
\begin{itemize}
\item[(a)] If each latent variable has exactly two pure observed children (that is, $\G=(\I_K;~ \I_K)^\top$ up to row permutations) and $H=2$, then there exist certain model parameters $\nnu$ and $\bo\Theta$ in the parameter space that are not identifiable.
\item[(b)]
If each latent variable has at least three pure observed children (i.e., $\G=(\I_K;~ \I_K;~ \I_K;~ \G_0^\top)^\top$ up to potential row permutations), then all the parameters $\nnu$ and $\bo\Theta$ in the parameter space are always identifiable up to permutations of the latent variables.
\end{itemize}
\end{proposition}

Proposition \ref{prop-cont}(a) shows that the sufficient condition for identifying the graph does not always suffice for identifying all the parameters in the model. Moreover, Proposition \ref{prop-cont}(b) shows that as long as the number of pure connections of each latent variable reaches three instead of two, then the identifiability of all parameters $\nnu$ and $\bo\Theta$ is guaranteed.
Along with Proposition \ref{prop-cont}, the main result Theorem \ref{thm-main} reveals an interesting phenomenon that identifying the discrete graph can require strictly weaker conditions than identifying the generative mechanism of the whole model.
The proof of Proposition \ref{prop-cont}(b) uses a similar argument as \cite{allman2009} by exploiting the classical Kruskal's theorem \citep{kruskal1977three} on unique three-way tensor decompositions. Identifiability proofs of this kind are not constructive, unlike our tensor unfolding-based technique.

\section{Discussion}\label{sec-disc}
In this paper, we have used a novel constructive proof based on unfolding tensors to identify the key loading graph in discrete latent bipartite graphical models. 
One notable feature of our identifiability result is that it holds regardless of whether the graph is directed or undirected, whether the latent variables are dependent or independent, and whether there are nonlinear transformations.
Our proof technique has an algebraic statistical nature, and it may be possible to generalize it to study other graphical models with discrete latent variables in the future.

The new identifiability condition requiring each latent variable to connect to at least two pure observed variables is practically interpretable.
In applications such as medical diagnosis or educational cognitive diagnosis, this condition can be interpreted as for each latent disease or latent skill, there exist at least two symptoms or test questions that only depend on this latent trait but not others.
In machine learning, our result can also provide some justification for adopting sparse loading structures in generative models such as RBMs.
Notably, we do not restrict the number of parents/neighbors the remaining non-pure observed variables can have, and they can have arbitrarily dense loadings on the latent layer. This ensures that our identifiable bipartite graphical models indeed define a very flexible model class beyond the tree-based models.

The identification result in this work illuminates a condition for \emph{population identifiability} to hold.
An important future direction is to develop a practical algorithm that takes an empirical tensor (instead of a population tensor considered here) based on a finite sample as input and outputs an estimator of the bipartite graph.
Further down the road, it would also be interesting to investigate the finite sample guarantees for such estimation algorithms.  
To address these questions, it might be useful to consider the nuclear norm relaxation of the matrix rank constraints or other spectral methods as considered previously for latent tree models \citep[e.g.][]{eriksson2005tree, ishteva2013unfolding}. 
We leave these directions for future research.

\bibliographystyle{apalike}
\bibliography{ref}

\clearpage
{\centering\section*{Supplementary Material}}
\addcontentsline{toc}{section}{Supplementary Material} 

\begin{appendix}
In Section A, we provide the proofs of the main theoretical results, including Theorem 1 and Propositions 1 and 2.
In Section B, we provide additional proofs of supporting lemmas and corollaries.
\section{Proofs of the Main Results}

We first state three useful lemmas.
As mentioned in the main text, for technical convenience we will slightly abuse notations when writing the Kronecker product and Khatri-Rao product, in that columns and rows in those products could be subject to permutations. Any of such permutations will not impact the rank properties of an matrix, so they do not affect any proof argument in this paper.

\begin{lemma}\label{lem-redun-latent}
Suppose $M\subseteq[p]$ is an index set satisfying that when $j$ ranges in $M$, the collection of of the $M$ sets $\con(j)$ are mutually disjoint sets. Then, we have
\begin{align}\label{eq-single-kron}
\MP(\Y_{M}\mid \A_{\con(M)}) = \bigotimes_{j\in M} \MP(Y_j\mid A_{\con(j)}).
\end{align}
Now consider an arbitrary set $M\subseteq [p]$. 
%Suppose  $S\subsetneqq[K]$. 
For any set ${S}$ such that $\con(M) \subseteq {S} \subseteq [K]$ with $\con(M) \neq {S}$, we have
\begin{align}\label{eq-lem-redun}
    \MP(\Y_{M}\mid \A_{{S}}) = \one_{H^{|{S}|-|\con(M)|}}^\top \bigotimes \MP(\Y_{M}\mid \A_{\con(M)}).
\end{align}
Moreover, for any $M_1,M_2\subseteq[p]$ and $M_1\cap M_2=\varnothing$ and any $S\supseteq \con(M_1) \cup \con(M_2)$, we have 
\begin{equation}\label{eq-property-kr}
\MP(\Y_{M_1\cup M_2}\mid \A_S) = \MP(\Y_{M_1}\mid \A_S) \bigodot \MP(\Y_{M_2}\mid \A_S).
\end{equation}
\end{lemma}

\begin{lemma}\label{lem-subset}
Suppose $S_1, S_2\subseteq [J]$ are two disjoint sets that satisfy $\con(S_1) \subseteq \con(S_2)$, $\con(S_1) \neq \con(S_2)$
%$\con(S_2)\setminus \con(S_1)$ is a singleton set, 
and $\MP(\Y_{S_1} \mid \A_{\con(S_1)})$ has full column rank $H^{|\con(S_1)|}$.
Then 
\begin{align}\label{eq-lem-large}
\rank(\MP(\Y_{S_1\cup S_2}\mid \A_{\con(S_2)}))
> H^{|\con(S_1)|}.
\end{align}

\end{lemma}

\begin{lemma}\label{lem-j1j2}
%\label{as-nodina}
For $j_1\neq j_2\in[J]$ with $|\con(j_1) \cup \con(j_2)| > 1$, the following holds for generic parameters:
%$\con(j_1)=\con(j_2)$ and $|\con(j_1)|>1$, 
% the following holds generically: $V^2\times H^{|\con(j_1) \cup \con(j_2)|}$ matrix 
$$\rank(\MP(Y_{\{j_1,j_2\}}\mid \A_{\con(j_1)\cup\con(j_2)})) > H.$$
%has rank larger than $H$ .
\end{lemma}

%%%
\begin{lemma}[Lemma 3.3 in \cite{stegeman2007kruskal}]\label{lem-khatri}
Consider two matrices $\A$ of size $I\times R$ and $\B$ of size $J\times R$.
\begin{itemize}
\item[(a)] If $\krank(\A)=0$ or $\krank(B)=0$, then $\krank(\A \bigodot \B) = 0$.
\item[(b)] If $\krank(\A)\geq 1$ and $\krank(B)\geq 1$, then 
$$\krank(\A \bigodot \B) \geq \min(\krank(\A) + \krank(\B) - 1, ~ R).$$
\end{itemize}
\end{lemma}

%%%%%%%%%%%%%%%%%%%%%%%%%%%%%%%%%%%%%%%%%%%%
\subsection{Proof of Proposition 1}

The size of the unfolding $[\mt]_{\{j_1,j_2\},\; :}$ is $V^2 \times V^{J-2}$.
If $\con(j_1) = \con(j_2) = \{k\}$ (that is, $g_{j_1,k} = g_{j_2,k} = 1$ and $g_{j_1,m} = g_{j_2,m} = 0$ for all $m\in[K]\setminus\{k\}$), then we can write $[\mt]_{\{j_1,j_2\},\; :}$ as 
\begin{align*}
    [\mt]_{\{j_1,j_2\},\; :} = 
    \underbrace{\MP(\Y_{\{j_1,j_2\}}\mid A_k)}_{V^2\times H} \cdot \underbrace{\MP(A_k,\; \A)}_{H \times H^K} \cdot \underbrace{\MP(\Y_{-\{j_1,j_2\}} \mid \A)^\top}_{H^K \times V^{J-2}}.
\end{align*}
Due to the above matrix factorization of the unfolded tensor, $\rank([\mt]_{\{j_1,j_2\},\; :}) \leq H$ must hold, which proves the ``if'' part of Proposition 1. Next consider the ``only if'' part.
Suppose $\con(j_1) = \con(j_2) = \{k\}$ does not hold for any $k\in [K]$, then $\con(j_1) \cup \con(j_2)$ is not a singleton set, so $|\con(j_1) \cup \con(j_2)| \geq 2$. 
In this case,
\begin{align}\label{eq-j1j2-nonpure}
    [\mt]_{\{j_1,j_2\},\; :} 
    &= 
    \underbrace{\MP(\Y_{\{j_1,j_2\}}\mid \A_{\con(j_1) \cup \con(j_2)})}_{\P_1:~ V^2\times H^{|\con(j_1) \cup \con(j_2)|}} \cdot \underbrace{\MP(\A_{\con(j_1) \cup \con(j_2)},\; \A)}_{\P_2: ~H^{\con(j_1) \cup \con(j_2)} \times H^K} \cdot \underbrace{\MP(\Y_{[J]\setminus\{j_1,j_2\}} \mid \A)^\top}_{\P_3: ~H^K \times V^{J-2}}\\ \notag
    &= \P_1 \cdot \P_2 \cdot \P_3
\end{align}
We next examine the rank of each of the three matrices $\P_1$, $\P_2$, and $\P_3$ on the right hand side of the above expression.

%%%
\emph{First, consider} $\P_1 = \MP(\Y_{\{j_1,j_2\}}\mid \A_{\con(j_1) \cup \con(j_2)})$, which is a $V^2 \times H^{|\con(j_1) \cup \con(j_2)|}$ matrix.
We apply Lemma \ref{lem-j1j2} to obtain $\rank(\P_1) > H$ holds for generic parameters.

\emph{Second, consider} $\P_2= \MP(\A_{\con(j_1) \cup \con(j_2)},\; \A)$. For generic distribution of the latent variables $\A$, the matrix $\P_2$ has rank $H^{|\con(j_1) \cup \con(j_2)|} \geq H^2 > H$. That is, $\P_2$ in \eqref{eq-j1j2-nonpure} generically has full rank $H^2$, which is greater than $H$.

\emph{Third, consider} $\P_3$.
we have an important observation that when $\G$ contains two disjoint copies of $\I_K$ and $|\con(j_1) \cup \con(j_2)| \geq 2$, the $(J-2) \times K$ submatrix $\G_{[J]\setminus\{j_1,j_2\},:}$ contains at least one submatrix of $\I_K$ after some row permutation.
To see this, note that if $|\con(j_1) \cup \con(j_2)| \geq 2$, then \emph{either} $\con(j_1)$ and $\con(j_2)$ are both singleton sets with $\con(j_1)\neq\con(j_2)$, \emph{or} at least one of $\con(j_1)$ and $\con(j_2)$ is not a singleton set. In either of these two cases, there exist a set $S \subseteq [J]\setminus\{j_1,j_2\}$ such that $|S|=K$ and $\G_{S,:} = \I_K$.
Applying Lemma \ref{lem-redun-latent} gives
$$
\MP(\Y_{S} \mid \A) = \bigotimes_{j\in S} \MP(Y_j\mid A_{\con(j)}),
$$
so $\MP(\Y_{S} \mid \A)$ generically has full column rank $H^K$.
Since $S \subseteq [J]\setminus\{j_1,j_2\}$, we can write
$$
\MP(\Y_{[J]\setminus\{j_1,j_2\}} \mid \A) = \MP(\Y_{S} \mid \A) \bigodot \MP(\Y_{[J]\setminus(S\cup \{j_1,j_2\})} \mid \A).
$$
% We next restate a useful lemma from \cite{stegeman2007kruskal}.
%%%
% \begin{lemma}[Lemma 3.3 in \cite{stegeman2007kruskal}]\label{lem-khatri}
% Consider two matrices $\A$ of size $I\times R$ and $\B$ of size $J\times R$.
% \begin{itemize}
% \item[(a)] If $\krank(\A)=0$ or $\krank(B)=0$, then $\krank(\A \bigodot \B) = 0$.
% \item[(b)] If $\krank(\A)\geq 1$ and $\krank(B)\geq 1$, then 
% $$\krank(\A \bigodot \B) \geq \min(\krank(\A) + \krank(\B) - 1, ~ R).$$
% \end{itemize}
% \end{lemma}
Next, we will use Lemma \ref{lem-khatri} to proceed with the proof.
One implication of Lemma \ref{lem-khatri}(b) is that for two matrices $\A$ and $\B$ with the same number of columns, if $\A$ has full column rank and $\B$ does not contain any zero columns, then $\krank(\A \bigodot \B) \geq \min(R+1-1,~ R) = R$, so $\A\bigodot\B$ has full column rank.
Note that $\MP(\Y_{S} \mid \A)$ has full column rank and $\MP(\Y_{[J]\setminus(S\cup \{j_1,j_2\})} \mid \A)$ does not contain any zero column vectors since each column of it is a conditional probability vector. So by Lemma \ref{lem-khatri}, we obtain that $\MP(\Y_{[J]\setminus\{j_1,j_2\}} \mid \A)$ has full column rank $H^K$ generically. So, $\P_3$ in \eqref{eq-j1j2-nonpure} has full row rank $H^K$ generically.

Going back to \eqref{eq-j1j2-nonpure}, since both $\P_2$ and $\P_3$ have full rank, the rank of $[\mathcal T]_{\{j_1,j_2\}, :}$ equals the rank of $\P_1$, which is greater than $H$. Now we have proved the ``only if'' part of the proposition that $\rank([\mathcal T]_{\{j_1,j_2\}, :}) > H$ generically if $\con(j_1)\cup \con(j_2)$ is not a singleton set.
The proof of Proposition 1 is complete.
\qedb
% \end{proof}

\subsection{Proof of Proposition 2}
Define notations $B_j$ and $C$:
$$
B_j = ([K]\setminus\{k\}) \cup \{j\},\quad
C = \{K+1,\ldots,2K\}.
$$
We first prove the ``only if'' part of Proposition 2.
It suffices to show that if $k\not\in\con(j)$, then $\rank([\mathcal T]_{B_j, C}) \leq H^{K-1}$.
Suppose $k\not\in\con(j)$, then $\con(j) \subseteq \{1,\ldots,k-1,k+1,\ldots,K\}$ and hence $B_j = \{1,\ldots,k-1,k+1,\ldots,K\}$.
In this case, we have the following decomposition,
\begin{align*}
    [\mathcal T]_{B_j, C} = \underbrace{\MP(\Y_{B_j}\mid \A_{[K]\setminus\{k\}})}_{V^K\times H^{K-1}} \cdot \underbrace{\MP(\A_{[K]\setminus\{k\}},~ \A_{[K]})}_{H^{K-1}\times H^K} \cdot \underbrace{\MP(\Y_{C}\mid \A_{[K]})^\top}_{H^K\times V^K}.
\end{align*}
Since the first matrix factor above has $H^{K-1}$ columns, the above decomposition clearly shows that $\rank([\mathcal T]_{B_j, C}) \leq H^{K-1}$.
%%%

We next prove the ``if'' part of Proposition 2.
Since $B_j = ([K]\setminus\{k\})\cup \{j\}$ and $k\in\con(j)$, we can rewrite $B_j = \big\{[K]\setminus\con(j)\big\} \bigcup \big\{(\con(j)\setminus\{k\})\cup \{j\})\big\}$. So we have
\begin{align}\label{eq-bjc}
    [\mathcal T]_{B_j, C} 
    &= \underbrace{\MP(\Y_{B_j}\mid \A_{1:K})}_{\text{denoted as }\P_1} \cdot \underbrace{\diag(\MP(\A_{1:K}))}_{\text{denoted as }\P_2} \cdot \underbrace{\MP(\Y_{C}\mid \A_{1:K})^\top}_{\text{denoted as }\P_3},\\ \notag
    \text{where }~
    \P_1 
    &=
     \MP(\Y_{[K]\setminus\con(j)}\mid \A_{[K]\setminus\con(j)}) \bigotimes 
    \MP(\Y_{(\con(j)\setminus\{k\})\cup \{j\})}\mid \A_{\con(j)})\\ \notag
    &=
    \underbrace{\Big(\bigotimes_{m\in[K]\setminus\con(j)} \MP(Y_m\mid A_m) \Big)}_{\text{denoted as }\P_{1,1}:~ V^{K-|\con(j)|} \times H^{K-|\con(j)|}} \bigotimes 
    \underbrace{\MP(\Y_{(\con(j)\setminus\{k\})\cup \{j\})}\mid \A_{\con(j)})}_{\text{denoted as }\P_{1,2}:~ V^{|\con(j)|} \times H^{|\con(j)|}}.
\end{align}
First, the matrix factor $\P_{1,1}$ is a Kronecker product of $K-|\con(j)|$ matrices, each of which has size $V\times H$ and has full column rank $H$ for generic parameters. So, $\P_{1,1}$ has full rank $H^{K-|\con(j)|}$ generically.
Since $\rank(\P_1) = \rank(\P_{1,1}) \cdot \rank(\P_{1,2})$, we next will prove 
\begin{align}\label{eq-rank12}
    \rank(\P_{1,2}) > H^{|\con(j)|-1}
\end{align} 
in order to show $\rank(\P_1) > H^{K-|\con(j)|} \cdot H^{|\con(j)|-1} = H^{K-1}$ in Proposition 2.
Denote by $\G_{\con(j)\cup \{j\},\; \con(j)}$ the submatrix of $\G$ with row indices ranging in $\con(j)\cup \{j\}$ and column indices ranging in $\con(j)$. Without loss of generality, after some column permutation $\G_{\con(j)\cup \{j\},\; \con(j)}$ can be written as follows and we denote it as $\widetilde\G$ for convenience:
\begin{align*}
    \widetilde\G=
    \G_{\con(j)\cup \{j\},\; \con(j)} = 
    \begin{pmatrix}
        0 & 1 & 0 & \cdots & 0 \\
        0 & 0 & 1 & \cdots & 0 \\
        \vdots & \vdots & \vdots & \ddots & \vdots \\
        0 & 0 & 0 & \cdots & 1 \\
        1 & 1 & 1 & \cdots & 1 
    \end{pmatrix};
\end{align*}
in other words, after some column permutation we can always write $\G_{\con(j)\cup \{j\},\; \con(j)}$ such that its first column corresponds to the latent variable $A_k$ under our current consideration, and its remaining columns correspond to all the other parent latent variables of $Y_j$.
Write $\P_{1,2}$ as
$$\P_{1,2} = \MP(\Y_{\con(j)\setminus\{k\}} \mid \A_{\con(j)}) \bigodot \MP(\Y_j \mid \A_{\con(j)}),$$
where the second matrix factor above $\MP(\Y_j \mid \A_{\con(j)})$ has size $V\times H^{|\con(j)|}$.
%%%%%
We next use Lemma \ref{lem-j1j2} to lower bound the rank of $\P_{1,2}$. Recall that $\G_{1:K,:} = \I_K$.
Define $S_1 = \con(j) \setminus \{k\}$ and $S_2 = \{j\}$, then $\con(S_1) = \con(j) \setminus \{k\}$ and $\con(S_2) = \con(j)$. The sets $S_1$ and $S_2$ satisfy that $\con(S_1) \subseteq \con(S_2)$ and $\con(S_2) \setminus \con(S_1)$ is a singleton set $\{k\}$; additionally, 
$$\MP(\Y_{S_1} \mid \A_{\con(S_1)}) = \bigotimes_{\ell\in \con(j) \setminus \{k\}} \MP(Y_\ell\mid A_{\ell})$$
has full column rank $H^{|\con(S_1)|}$.
Therefore, we can use Lemma \ref{lem-j1j2} to obtain 
$
\rank(\P_{1,2}) > H^{|\con(S_1)|}.
$
According to the argument right after \eqref{eq-rank12}, it holds that
$$
\rank(\P_1) > H^{K-1}.
$$
We now go back to the decomposition in \eqref{eq-bjc} that $[\mathcal T]_{B_j, C} = \P_1 \P_2 \P_3$. The rank of the $H^K \times H^K$ diagonal matrix $\P_2 = \diag(\MP(\A_{1:K}))$ is equal to $H^K$ for generic parameters.
The rank of the $H^K \times V^K$ matrix $\P_3$ is also equal to $H^K$ for generic parameters, because $\P_3 = \bigotimes_{j=K+1}^{2K} \MP(Y_j\mid \A_{j-K})^\top$. 
% Matrix $[\mathcal T]_{B_j, C}$ has size $V^K\times V^K$.
% Since $V\geq H$, so 
Both $\P_2$ and $\P_3$ have full rank generically by Assumption 2.
Therefore, $\rank([\mathcal T]_{B_j, C}) = \rank(\P_1 \P_2 \P_3) = \rank(\P_1) > H^{K-1}$ generically. This proves the ``if'' part of Proposition 2.
\qedb

%%%%%%%%%%%%%%%%%

\subsection{Proof of Theorem 1}
Let $\mathcal T = (t_{i_1,\ldots,i_J})_{i_1,\ldots,i_J \in\{0,1,\ldots,V-1\}}$ denote the $J$-mode probability tensor with $V^J$ entries that characterizes the marginal joint distribution of the $J$ observed variables. Proposition 1 has the following implication. If we exhaustively search for all pairs of observed variables when $\{j_1,j_2\}\subseteq [p]$ with $j_1\neq j_2$, then we can exactly obtain $K$ disjoint sets of variable indices $S_1,\ldots, S_{K}\subseteq [p]$ with the following property: For any $k\in[K]$, it holds that
\begin{align*}
    \rank([\mt]_{\{j_1,j_2\},\; :}) \leq H,\quad \text{for any } j_1, ~ j_2\in S_k.
\end{align*}
We next show that the above argument can also uniquely identify the number of latent variables $K$. This is true because $K$ is equal to the number of all the nonoverlapping sets $S_1, S_2,\ldots \subseteq [J]$ such that for any $j_1,j_2\in S_k$ for each $k$, it holds that $\rank([\mathcal T]_{\{j_1,j_2\},:}) \leq H$.
After exhaustively inspecting $\rank([\mathcal T]_{\{j_1,j_2\},:})$ for all pairs of indices $j_1,j_2\in[J]$ for the observed variables, the number of such nonoverlapping sets is uniquely determined and $K$ is uniquely determined in a constructive way.
Thus far, we have identified all the row vectors indexed by $j\in\bigcup_{k=1}^K S_k$. Each of these row vectors is a standard basis vector.
Now using the assumption that the true graphical matrix $\G$ satisfies that each column contains at least two entries of ``1'', we can swap the orders of the already identified row vectors in $\G$ so that $\G$ takes the form of $\G=(\I_K; \I_K, \G^{*\top})^\top$.
Note that in this setting, we have $\con(k) = k$ for all $k\in[K]$.

Further, for any $j\in[J]$ with $|\con(j)| > 1$, in order to determine whether $k\in\con(j)$ for each $k\in[K]$, we consider the following unfolding 
\begin{align}\notag
%\label{eq-setab}
[\mathcal T]_{B_{j,k}, C},~\text{ where }
B_{j,k} := ([K]\setminus\{k\}) \cup \{j\},\quad
C := \{K+1,\ldots,2K\}.
\end{align}
Here both $B_{j,k}$ and $C$ are non-overlapping sets with cardinality $K$, so $[\mathcal T]_{B_{j,k}, C}$ is a $H^K \times H^{K}$ matrix.
If $|\con(j)|\geq 2$, Proposition 2 implies that $\rank([\mathcal T]_{B_{j,k}, C}) > H^{K-1}$ holds if and only if $k\in\con(j)$. Therefore, examining whether the rank of $\rank([\mathcal T]_{B_{j,k}, C})$ exceeds $H^{K-1}$ exhaustively for all $j$ with $|\con(j)|\geq 2$ and all $k\in[K]$ will uniquely identify the entries in the $\G$ matrix in the rows indexed by those $j$ with $|\con(j)|\geq 2$.
This constructively proves the identifiability of the row vectors in $\G$ indexed by those $j$ for which $|\con(j)| \geq 2$. 
Combined with the fact established in the previous paragraph that all row vectors in $\G$ indexed by those $j$ for which $|\con(j)|=1$ are identified, we complete the proof that the entire $\G$ matrix is identified.
The proof of Theorem 1 is complete.
\qedb

\section{Proof of Supporting Results}
\subsection{Proof of Lemma \ref{lem-redun-latent}}
We first prove \eqref{eq-single-kron}.
By the conditional independence of the observed variables given the latent parents, we have
\begin{align*}
    \MP(\Y_M\mid \A_{\con(M)}) 
    = \prod_{j=1}^M \MP(Y_j\mid A_{\con(M)})
    = \prod_{j=1}^M \MP(Y_j\mid A_{\con(j)}),
\end{align*}
which implies for any observed pattern $\bo y_M\in\{0,\ldots,V-1\}^{|M|}$ and any latent pattern $\bo a_{\con(M)} \in \{0,\ldots,H-1\}^{|\con(M)|}$, it holds that
\begin{align*}
    \MP(\Y_M = \bo y_M\mid \A_{\con(M)} = \bo a_{\con(M)})
    = \prod_{j=1}^M \MP(Y_j = y_j\mid A_{\con(j)}=a_{\con(j)}).
\end{align*}
Since $\con(j)$ for $j\in M$ are disjoint singleton sets, by following the definition of the Kronecker product of matrices, the above factorization implies $\MP(\Y_{M}\mid \A_{\con(M)}) = \bigotimes_{j\in M} \MP(Y_j\mid A_{\con(j)})$ and proves \eqref{eq-single-kron}.

We next prove \eqref{eq-lem-redun}.
% Introduce notation $\widetilde{S} = \con(M)$.
First, note that $\MP(\Y_{M}\mid \A_{{S}})$ has size $V^{|M|} \times H^{|{S}|}$, while $\one_{H^{|{S}|-|\con(M)|}}^\top$ has size $1\times H^{|S|-|\con(M)|}$ and $\MP(\Y_{M}\mid \A_{\con(M)})$ has size $V^{|M|} \times H^{|\con(M)|}$. So, the sizes of matrices on the left hand side and the right hand side of \eqref{eq-lem-redun} match each other. 
%%%
Further, since $\con(M) \subsetneqq S$, we know that the conditional distribution of $\Y_M$ given $\A_{S}$ only depends on those latent variables belonging to the index set $\con(M)$; in other words, 
$$\MP(\Y_M\mid \A_{S}) = \MP(\Y_{M}\mid \A_{\con(M)}, \A_{S\setminus \con(M)}) = \MP(\Y_{M}\mid \A_{\con(M)}),
$$
where $\A_{\con(M)}$ and $\A_{S\setminus \con(M)}$ are a subvectors of $\A_{S}$.
%%%
The above fact implies that, although $\MP(\Y_{M}\mid \A_{S})$ has $H^{|S|}$ columns indexed by all possible configurations of $\A_{S}\in\{0,\ldots,H-1\}^{|S|}$, actually $\MP(\Y_{M}\mid \A_{S})$ contains only $H^{|\con(M)|} (<H^{|S|})$ unique columns. Each unique column is indexed by a unique pattern the subvector $\A_{\con(M)}$ takes in $\{0,\ldots,H-1\}^{|\con(M)|}$. So, the matrix $\MP(\Y_M\mid \A_{S})$ horizontally stacks $H^{|S|-|\con(M)|}$ repetitive blocks, each of which equals $\MP(\Y_{M}\mid \A_{\con(M)})$. So, we have
\begin{align*}
    \MP(\Y_{M}\mid \A_{S})
    &=
    \underbrace{
    \begin{pmatrix}
    \MP(\Y_{M}\mid \A_{\con(M)}), & \cdots, & \MP(\Y_{M}\mid \A_{\con(M)}) 
    \end{pmatrix}
    }_{H^{|S|-|\con(M)|}\text{ copies}}   
    \\
    &=
    \one_{H^{|S|-|\con(M)|}}^\top \bigotimes \MP(\Y_{M}\mid \A_S).
\end{align*}

We next prove \eqref{eq-property-kr}. Since $S\supseteq \con(M_1) \cup \con(M_2)$ holds, $\Y_{M_1}$ and $\Y_{M_2}$ are conditionally independent given $\A_{S}$. So, we can use Lemma 12 in \cite{allman2009} to direct obtain the conclusion that $\MP(\Y_{M_1\cup M_2}\mid \A_S) = \MP(\Y_{M_1}\mid \A_S) \bigodot \MP(\Y_{M_2}\mid \A_S)$.
The proof of Lemma \ref{lem-redun-latent} is complete.
\qedb

\subsection{Proof of Lemma \ref{lem-subset}}
Since $\con(S_1) \subsetneqq \con(S_2)$, we use Lemma \ref{lem-redun-latent} to write
\begin{align*}
    \MP(Y_{S_1} \mid \A_{\con(S_2)}) = \one_{H^{|\con(S_2)\setminus \con(S_1)|}}^\top \bigotimes \MP(Y_{S_1} \mid \A_{\con(S_2)}),
\end{align*}
and
\begin{align}\label{eq-s1s2}
    \MP(\Y_{S_1\cup S_2}\mid \A_{\con(S_2)}) = \Big(\one_{H^{|\con(S_2)\setminus \con(S_1)|}}^\top \bigotimes \MP(Y_{S_1} \mid \A_{\con(S_1)})\Big)
    \bigodot \MP(\Y_{S_2}\mid \A_{\con(S_2)}).
\end{align}
Define notation
\begin{align*}
\F &=  \MP(Y_{S_1} \mid \A_{\con(S_1)});\\
\E_{\aaa_\ell} 
&= \MP(\Y_{S_2}\mid \A_{\con(S_1)}, \A_{\con(S_2)\setminus \con(S_1)} = \aaa_\ell),\quad \aaa_\ell\in\{0,1,\ldots,H-1\}^{|\con(S_2)\setminus \con(S_1)|}.
\end{align*}
For each fixed pattern $\aaa_\ell\in\{0,1,\ldots,H-1\}^{|\con(S_2)\setminus \con(S_1)|}$, the $\E_{\aaa_\ell}$ is 
%%%%
a $V^{|S_2|} \times H^{|\con(S_1)|}$ matrix with rows indexed by all possible configurations of $Y_{S_2}$ and columns by all possible configurations of $\A_{\con(S_1)}$.
By the condition in the lemma, matrix $\F$ has full column rank $H^{|\con(S_1)|}$. 
For notational simplicity, denote this number $H^{|\con(S_1)|} =: R$.
After possibly rearranging the columns of $\one_{H^{\con(S_2)\setminus \con(S_1)}}^\top \bigotimes \MP(Y_{S_1} \mid \A_{\con(S_1)})$ and $\MP(\Y_{S_2}\mid \A_{\con(S_2)})$ via a common permutation, we can write 
\begin{align*}
    \one_{H^{\con(S_2)\setminus \con(S_1)}}^\top \bigotimes \MP(Y_{S_1} \mid \A_{\con(S_1)}) 
    &= 
    (\underbrace{\F, \F, \ldots, \F}_{R=H^{\con(S_2)\setminus \con(S_1)} \text{ copies}});\\
    \MP(\Y_{S_2}\mid \A_{\con(S_2)})
    &= (\E_{\aaa_1}, \E_{\aaa_2}, \ldots, \E_{\aaa_R}).
\end{align*}
Therefore, the matrix $\MP(\Y_{S_1\cup S_2}\mid \A_{\con(S_2)})$ in \eqref{eq-s1s2} can be written as 
\begin{align}\label{eq-s1s2-ef}
    \MP(\Y_{S_1\cup S_2}\mid \A_{\con(S_2)})
    = 
    \Big(\F \bigodot \E_{\aaa_1},~ \F \bigodot \E_0,~ \ldots,~ \F \bigodot \E_{\aaa_R} \Big).
\end{align}
Next, we
will prove that the matrix $\MP(\Y_{S_1\cup S_2}\mid \A_{\con(S_2)})$ has rank at least $R + 1$ under Assumption 2.
Denote the columns of $\F$ by $\F_1,\ldots,\F_R$, and the columns of $\E_{\aaa_h}$ by $\E_{\aaa_h,1},\ldots,\E_{\aaa_h,R}$.
Note that by definition, each $\E_{\aaa_h,m}$ is a vector describing the conditional distribution of $Y_{S_2}$ given a specific configuration of $\A_{\con(S_2)}$ when $\A_{\con(S_2)\setminus\con(S_1)}$ is fixed to be $\aaa_h$.
Since $\con(S_1) \subsetneqq \con(S_2)$, there exists some $k\in \con(S_1)$ but $k\not\in\con(S_2)$.
By Assumption 2, 
there exist $H$ different vectors $\bo\beta_1,\ldots,\bo\beta_H \in \{\aaa_1,\ldots,\aaa_R\}$ where $\bo\beta_1,\ldots,\bo\beta_H$ differ only in the entry corresponding to $A_k=0,1,\ldots, H-1$, in which
%the vectors the $H$ vectors $\E_{\bo\beta_1}, \ldots, \E_{\bo\beta_H}$ are linearly independent.
% So there must 
there
exist $h_1\neq h_2 \in\{0,1,\ldots,H-1\}$ and some $r\in\{1,\ldots,R\}$ such that 
\begin{equation}\label{eq-h1h2}
    \E_{\bo\beta_{h_1},r} \neq \E_{\bo\beta_{h_2},r}.
\end{equation}
Suppose that there exists a vector $\bo x = (x_1,\ldots,x_{R}, x_{R+1})^\top$ such that
\begin{align*}
    x_1 \F_1 \bigotimes \E_{\bo\beta_{h_1},1} + \cdots +  x_R \F_R \bigotimes \E_{\bo\beta_{h_1},R} + x_{R+1} \F_{r} \bigotimes \E_{\bo\beta_{h_2}, r} &= \zero.
\end{align*}
Rearranging the terms above gives
\begin{align}\label{eq-fe}
    \sum_{m\in[R]\setminus\{r\}} x_m \F_m \bigotimes \E_{\bo\beta_{h_1},m} + \F_r \bigotimes (x_r \E_{\bo\beta_{h_1},r} + x_{R+1}\E_{\bo\beta_{h_2},r}) = \zero.
\end{align}
Define a matrix $\widetilde \E$ whose $r$-th column vector is $x_r \E_{\bo\beta_{h_1},r} + x_{R+1}\E_{\bo\beta_{h_2},r}$, and for any $m\in[R]$ and $m\neq r$, the $m$-th column vector of $\widetilde\E$ is $\E_{\bo\beta_{h_1},m}$. Define a $R$-dimensional vector $\widetilde{\bo y}$ with $y_r = 1$ and $y_m=x_m$ for any $m\in[R]$ and $m\neq r$.
Then \eqref{eq-fe} can be rewritten as
\begin{align}\label{eq-fe-odot}
    (\F \bigodot \widetilde\E) \bo y = \zero.
\end{align}
Since $\bo y$ above is not a zero vector, we have that $\rank(\F \bigodot \widetilde\E) < R$.
We next will show that this necessarily implies $x_r \E_{\bo\beta_{h_1},r} + x_{R+1}\E_{\bo\beta_{h_2},r} = \zero$.

We introduce the definition of the Kruskal rank.
The Kruskal rank of a matrix is the largest integer $M$ such that every $M$ columns of this matrix are linearly independent. Denote the Kruskal rank of any matrix $\B$ by $\krank(\B)$. It is easy to see that if a matrix has full column rank, then its Kruskal rank equals its rank. Also, a matrix $\mathbf B$ contains a zero column vector if and only if $\krank(\mathbf B)=0$.

Now consider the matrix $\F \bigodot \widetilde\E$ in \eqref{eq-fe-odot}. Since $\F$ has full column rank $R$, then $\krank(\F) = \rank(\F) = R$. 
For any $m\in[R]\setminus\{r\}$, the $m$-th column vector of $\widetilde\E$ is a conditional probability vector for $Y_j$ given a fixed pattern of $\A_{\con(j)}$, therefore $\E_{\bo\beta_{h_1},m}\neq \zero$. So, if $x_r \E_{\bo\beta_{h_1},r} + x_{R+1}\E_{\bo\beta_{h_2},r} \neq \zero$, then $\widetilde\E$ does not have any zero column vectors and $\krank(\widetilde\E) \geq 1$. In this case, the following holds by Lemma \ref{lem-khatri}(b),
\begin{align*}
    \krank(\F \bigodot \widetilde\E) 
    &\geq \min(\krank(\F) + \krank(\widetilde\E) - 1, ~ R)\\
    &= \min(R + \krank(\widetilde\E) - 1, ~ R) \\
    &\geq \min(R + 1 - 1, ~ R) = R.
\end{align*}
Then, $\rank(\F \bigodot \widetilde\E) = \krank(\F \bigodot \widetilde\E) = R$, which contradicts with \eqref{eq-fe-odot} because the vector $\bo y$ there is not a zero vector with $y_r=1$ and hence $\F \bigodot \widetilde\E$ should not have full column rank.
This contradiction implies that 
$$
x_r \E_{\bo\beta_{h_1},r} + x_{R+1}\E_{\bo\beta_{h_2},r} = \zero.
$$
Now note that each of $\E_{\bo\beta_{h_1},r}$ and $\E_{\bo\beta_{h_2},r}$ is a conditional probability vector describing the distribution of $Y_{S_2}$ given a fixed pattern of $\A_{\con(S_2)}$. Therefore, $\one^\top\E_{\bo\beta_{h_1},r} = \one^\top\E_{\bo\beta_{h_1},r} = 0$, and we further have $x_r + x_{R+1} = 0$.
Plugging the above expression $x_r \E_{\bo\beta_{h_1},r} + x_{R+1}\E_{\bo\beta_{h_2},r} = \zero$ back into \eqref{eq-fe} further gives
\begin{align}\label{eq-fe-minus1}
    \sum_{m\in[R]\setminus\{r\}} x_m \F_m \bigotimes \E_{\bo\beta_{h_1},m} = \zero.
\end{align}
We next use Lemma \ref{lem-khatri} to show that $\{\F_m \bigotimes \E_{\bo\beta_{h_1},m}:~ m\in[R] \setminus\{r\}\}$ are $R-1$ linearly independent vectors. To show this, first note that \eqref{eq-fe-minus1} can be written as 
$$
(\F_{-r} \bigodot \widetilde \E_{-r}) \bo x_{-r} = \zero,
$$
where matrix $\F_{-r}$ contains all but the $r$-th column of matrix $\F$, matrix $\widetilde \E_{-r}$ contains all but the $r$-th column of matrix $\widetilde\E$, and vector $\bo x_{-r}$ contains all but the $r$-th entry of the vector $\bo x$.
Now, because $\krank(\F_{-r}) = \rank(\F_{-r}) = R-1$ and $\krank(\widetilde \E_{-r}) \geq 1$, Lemma \ref{lem-khatri}(b) implies that 
$$
\krank(\F_{-r} \bigodot \widetilde \E_{-r}) \geq \min(R-1 + \krank(\widetilde \E_{-r})-1, ~R-1) \geq R-1.
$$
Since $\F_{-r} \bigodot \widetilde \E_{-r}$ contains $R-1$ columns, the above inequality gives that $\rank(\F_{-r} \bigodot \widetilde \E_{-r}) = \krank(\F_{-r} \bigodot \widetilde \E_{-r}) = R-1$.
In other words, the $R-1$ vectors in $\{\F_m \bigotimes \E_{\bo\beta_{h_1},m}:~ m\in[R]\setminus\{r\}\}$ are linearly independent. Therefore, from \eqref{eq-fe-minus1} we have that $x_m = 0$ for all $m\in[R]\setminus\{r\}$.
Recall we aim to examine whether there exists a nonzero vector $\bo x = (x_1,\ldots,x_R,x_{R+1})$ such that \eqref{eq-fe} holds, and we have already shown that under \eqref{eq-fe}, the only two potentially nonzero entries in the vector $\bo x$ are $x_r$ and $x_{R+1}$. 
If $(x_r, x_{R+1}) \neq (0,0)$, then $x_r = - x_{R+1} \neq 0$, and it indicates 
$$
\E_{\bo\beta_{h_1},m} = \E_{\bo\beta_{h_2},m}.
$$
The above equality contradicts the earlier \eqref{eq-h1h2} that $\E_{\bo\beta_{h_1},m} \neq \E_{\bo\beta_{h_2},m}$ which follows from Assumption 2, therefore $x_r=x_{R+1}=0$ must hold. 
In summary, we have shown that the $R+1$ vectors $\F_1 \bigotimes \E_{\bo\beta_{h_1},1},~ \ldots,~ \F_R \bigotimes \E_{\bo\beta_{h_1},R},~ \F_{r} \bigotimes \E_{\bo\beta_{h_2}, r}$ are linearly independent, which means the matrix $\MP(Y_{S_1\cup S_2}\mid \A_{\con(S_2)})$ in \eqref{eq-s1s2-ef} contains at least $R+1$ linearly independent column vectors and $\rank(\MP(Y_{S_1\cup S_2}\mid \A_{\con(S_2)})) \geq R+1 = H^{|\con(S_1)|} + 1$.
The proof of Lemma \ref{lem-subset} is complete.

%%%%%%%%%%%%%%%%%%%
\subsection{Proof of Lemma \ref{lem-j1j2}}
Consider all possible cases: 
Case (a), $\con(j_1) \cap \con(j_2) = \varnothing$; 
Case (b), $\con(j_1) = \con(j_2)$ with $|\con(j_1)| = |\con(j_2)| \geq 2$;
Case (c), $\con(j_1) \subsetneqq \con(j_2)$ with $|\con(j_2)| \geq 2$;
and Case (d), $\con(j_1) \not\subseteq \con(j_2)$ and $\con(j_2) \not\subseteq \con(j_1)$ with $|\con(j_1)|\geq 2$,  $|\con(j_2)| \geq 2$.
Denote $\P = \MP(Y_{\{j_1,j_2\}}\mid \A_{\con(j_1)\cup\con(j_2)})$.

\noindent\textbf{Case (a):} Since $\con(j_1) \cap \con(j_2) = \varnothing$, we have 
\begin{align*}
    \P= \MP(Y_{\{j_1,j_2\}}\mid \A_{\con(j_1)\cup\con(j_2)}) = \MP(\Y_{j_1}\mid \A_{\con(j_1)}) \bigotimes \MP(\Y_{j_2}\mid \A_{\con(j_2)}),
\end{align*}
so $\rank(\P) = \rank(\MP(\Y_{j_1}\mid \A_{\con(j_1)})) \cdot \rank(\MP(\Y_{j_2}\mid \A_{\con(j_2)})) \geq H^2 > H$.

\bigskip
\noindent\textbf{Case (b):} We have $\con(j_1) \cup \con(j_2) = \con(j_1) = \con(j_2)$, so both $\MP(\Y_{j_1}\mid \A_{\con(j_1)}) = \MP(\Y_{j_1}\mid \A_{\con(j_1) \cup \con(j_2)})$ and $\MP(\Y_{j_2}\mid \A_{\con(j_2)}) = \MP(\Y_{j_2}\mid \A_{\con(j_1) \cup \con(j_2)})$ are generic conditional probability matrices without particular structures (i.e., without equality between certain columns). We can write
$$
\P=
\MP(\Y_{j_1}\mid \A_{\con(j_1) \cup \con(j_2)}) \bigodot \MP(\Y_{j_2}\mid \A_{\con(j_1) \cup \con(j_2)}).
$$
Now we apply Lemma 13 in \cite{allman2009} to obtain that $\rank(\P) = \min(V^2, H^{|\con(j_1) \cup \con(j_2)|}) = H^{|\con(j_1) \cup \con(j_2)|} > H$.

\bigskip
\noindent\textbf{Case (c):}
denote $L = H^{|\con(j_2)\setminus\con(j_1)|}$. We can apply Lemma \ref{lem-redun-latent} to write
\begin{align}\notag
\P 
&=
\MP(\Y_{j_1}\mid \A_{\con(j_2)}) \bigodot \MP(\Y_{j_2}\mid \A_{\con(j_2)})\\
\notag
&= \Big(\one_{H^{|\con(j_2)\setminus\con(j_1)|}}^\top \bigotimes \MP(\Y_{j_1}\mid \A_{\con(j_1)})\Big) \bigodot \MP(\Y_{j_2}\mid \A_{\con(j_2)}) \\
\label{eq-eal}
&= \Big(\MP(\Y_{j_1}\mid \A_{\con(j_1)}) \bigodot \E_{\aaa_1},~ \ldots,~ \MP(\Y_{j_1}\mid \A_{\con(j_1)}) \bigodot \E_{\aaa_L}\Big),
\end{align}
where each $\E_{\aaa_\ell}$ denotes a $V \times H^{|\con(j_1)|}$ submatrix of $\MP(\Y_{j_2}\mid \A_{\con(j_2)})$. Here, $\aaa_\ell\in\{0,1,\ldots,H-1\}^{|\con(j_2) \setminus \con(j_1)|}$ indexes a specific configuration of $\A_{\con(j_2) \setminus \con(j_1)}$, and the $H^{|\con(j_1)|}$ columns of $\E_{\aaa_\ell}$ are indexed by all possible configurations $\A_{\con(j_1)}$ can take.
Therefore, each matrix $\E_{\aaa_\ell}$ collections the conditional probability vectors of $Y_{j_2}$ when fixing $\A_{ \con(j_2) \setminus \con(j_1)} = \aaa_\ell$ and varying $\A_{\con(j_1)}$.
Note that $\E_{\aaa_\ell}$ is a generic matrix without equality constraints between certain columns. Therefore, we apply Lemma 13 in \cite{allman2009} to obtain that the following holds for generic parameters:
$$
\rank(\MP(\Y_{j_1}\mid \A_{\con(j_1)}) \bigodot \E_{\aaa_\ell}) = \min(V^2, H^{|\con(j_1)|}).
$$
Now, if $|\con(j_1)|>1$, then the above inequality gives $\rank(\MP(\Y_{j_1}\mid \A_{\con(j_1)}) \bigodot \E_{\aaa_\ell}) \geq H^2 > H$, which further implies $\rank(\P) \geq \rank(\MP(\Y_{j_1}\mid \A_{\con(j_1)}) > H$.
So, it remains to consider the case with $|\con(j_1)|=1$. In this case, $\MP(\Y_{j_1}\mid \A_{\con(j_1)})$ has full column rank $H$ by Assumption 2.
Define $S_1=\{j_1\}$ and $S_2=\{j_2\}$, then the conditions $S_1\subsetneqq S_2$ and $\MP(\Y_{S_1}\mid \A_{\con(S_1)})$ in Lemma \ref{lem-subset} are satisfied.
So Lemma \ref{lem-subset} gives that the rank of $\P$ is greater than $H^{|\con(S_1)|} = H$.
In summary, in Case (c), $\rank(\P) > H$.

\bigskip
\noindent\textbf{Case (d):} $\con(j_1)\subsetneqq\con(S_2)$ and $\con(j_2)\subsetneqq\con(S_1)$ imply that both $\con(j_1)$ and $\con(j_2)$ contain at least two elements.
First consider the $V^{2} \times H^{|\con(j_2)|}$ matrix
$$
\MP(\Y_{\{j_1,j_2\}}\mid \A_{\con(j_2)}) = \one_{H^{|\con(j_2) \setminus \con(j_1)|}}^\top \bigotimes \MP(\Y_{j_1}\mid\A_{\con(j_1) \cap \con(j_2)}) \bigodot \MP(\Y_{j_2}\mid\A_{\con(j_2)}),
$$
which takes the same form as \eqref{eq-eal} in case (c) and falls into the same setting considered in case (c). Therefore, we know that $\rank(\MP(\Y_{\{j_1,j_2\}}\mid \A_{\con(j_2)})) > H$ holds generically. Now, we claim that this matrix $\MP(\Y_{\{j_1,j_2\}}\mid \A_{\con(j_2)})$ has the following linear transformation relationship with the original matrix
$\P = \MP(\Y_{\{j_1,j_2\}}\mid \A_{\con(j_1) \cup \con(j_2)})$:
\begin{align}\label{eq-smatrix}
    \MP(\Y_{\{j_1,j_2\}}\mid \A_{\con(j_2)}) = \MP(\Y_{\{j_1,j_2\}}\mid \A_{\con(j_1) \cup \con(j_2)}) \cdot \mathbf S = \P\cdot \mathbf S,
\end{align}
where $\mathbf S$ has mutually orthogonal columns and hence full column rank.
To see that this claim is true, we only need to note that for any vectors $\bo y_{\{j_1,j_2\}} \in \{0,\ldots,V-1\}^{2}$ and $\aaa_{\con(j_2)}\in\{0,\ldots,H-1\}^{|\con(j_2)|}$, it holds that
\begin{align*}
&~ \MP(\Y_{\{j_1,j_2\}} = \bo y_{\{j_1,j_2\}} \mid \A_{\con(j_2)} = \aaa_{\con(j_2)}) \\
=&~ 
\sum_{\aaa_{\con(j_1)} \in\{0,\ldots,H-1\}^{|\con(j_1)|}}
\MP(\Y_{\{j_1,j_2\}} = \bo y_{\{j_1,j_2\}} \mid \A_{\con(j_2)} = \aaa_{\con(j_2)}, \A_{\con(j_1)} = \aaa_{\con(j_1)}).
\end{align*}
This equality means that to get the smaller conditional probability table $\MP(\Y_{\{j_1,j_2\}}\mid \A_{\con(j_2)})$, we can just sum up appropriate columns in the larger conditional probability table $\MP(\Y_{\{j_1,j_2\}}\mid \A_{\con(j_1) \cup \con(j_2)})$, where the $H^{|\con(j_1) \cup \con(j_2)|} \times H^{|\con(j_2)|}$ matrix $\mathbf S$ has binary entries that reflect this summation. Specifically, for any $\aaa_{\con(j_1)\cup\con(j_2)}$ and $\bo\beta_{\con(j_2)}$, the corresponding entry in $S$ is
$$
S_{\aaa_{\con(j_1)\cup\con(j_2)}, \bo\beta_{\con(j_2)}} = \mathbbm{1}(\aaa_{\con(j_1)\cup\con(j_2)}\text{ restricted to }\con(j_2)\text{ equals }
\bo\beta_{\con(j_2)}).
$$
Since each column in $\mathbf S$ is indexed by a different configuration $\bo\beta_{\con(j_2)}$, so the column vectors of $\mathbf S$ have mutually disjoint support and are orthogonal.
Now we have shown that the claim about  \eqref{eq-smatrix} is correct.
Therefore, \eqref{eq-smatrix} implies that 
$$\rank(\P) \geq \rank( \MP(\Y_{\{j_1,j_2\}}\mid \A_{\con(j_2)})) > H.$$
Now we have proved the conclusion of Lemma \ref{lem-j1j2} for all possible cases (a), (b), (c), and (d).
The proof is complete.

\subsection{Proof of Corollary 1}
% Consider $S \subseteq [J]\setminus\{j_1,j_2\}$. 
If $\con(j_1) = \con(j_2) = \{k\}$ for some $k\in[K]$, then the proof of Proposition 1 implies the following for any $S \subseteq [J]\setminus\{j_1,j_2\}$:
\begin{equation}\label{eq-sets}
\rank([\mathcal T]_{\{j_1, j_2\},\; S})
\leq  H,\quad \forall S \subseteq [J]\setminus\{j_1,j_2\} \text{ with } |S| \geq 2. 
\end{equation}

If $\con(j_1) = \con(j_2) = \{k\}$ does not hold, then
among all the possible $S \subseteq [J]\setminus\{j_1,j_2\}$ with $|S| \leq K$, there must exist some set $S$ such that $\rank(\unfold(\mathcal T(\{j_1, j_2\},\; S))) > H$.
In fact, we can just take $S$ to be a complete set of pure children of the $K$ latent variables. This is because as stated in the proof of Proposition 1, when $\con(j_1) = \con(j_2) = \{k\}$ does not hold, the matrix $\G_{[J]\setminus\{j_1,j_2\},:}$ contains at least one identity submatrix $\I_K$ after some row permutation;
we index this submatrix by $\G_{S,:}$ where $S \subseteq [J]\setminus\{j_1,j_2\}$ and $\G_{S,:} = \I_K$.
The proof of Proposition 1 gives that 
$
\rank([\mathcal T]_{\{j_1, j_2\},\; S}) > H
$
holds generically for this specific set $S$. The above inequality combined with \eqref{eq-sets} implies that
it suffices to consider subsets $S$ of cardinality $K$ and examine whether $[\mathcal T]_{\{j_1, j_2\},\; S} > H$ holds to determine the pure variable structure.
This proves Corollary 1.
% \end{proof}
\qedb

\subsection{Proof of Proposition 3}
To prove part (a) of the proposition, consider $H=2$ with binary latent variables and $\G=(\I_K;\I_K)^\top$. In this case, the latent bipartite graphical model reduces to the so-called BLESS model in \cite{gu2024blessing} with a star-forest graph between the observed layer and the latent layer. Therefore, we can use Theorem 3.3 in \cite{gu2024blessing} to obtain that if $A_k$ is independent of $\A_{[K]\setminus k}$ for some $k\in[K]$ (which can be translated to certain polynomial equation constraints on the proportion parameters $\nnu=(\nu_{\aaa};\aaa\in\{0,1\}^K)$), then the model parameters $(\bo\Theta,\nnu)$ are not identifiable. This proves the conclusion of part (a) of Proposition 3.

To prove part (b) of the proposition, we follow a similar argument as \cite{allman2009} by invoking the Kruskal's theorem \citep{kruskal1977three}. 
The proof is also similar to the proof of Proposition 3.4 in \cite{gu2024blessing} for the simpler BLESS model.
Suppose without loss of generality that {$\G=(\I_K,\;\I_K,\;\I_K,\;\G^{\star\top})^\top$}, where the submatrix $\G^\star$ can take an arbitrary form or even be empty. Suppose the alternative parameters $\ov\btheta$ and $\ov\nnu$ give rise to the same marginal distribution of the $J$ observed variables as the true parameters $\btheta$ and $\nnu$.
We group the first $K$ observed variables $Y_1,\ldots,Y_K$ into one categorical variable with $V^K$ categories and denote it by $Z_1$, then each of the $V^K$ possible configurations of the vector $\Y_{1:K}=(Y_1,\ldots,Y_K)$ corresponds to one category that $Z_1$ can take.
Similarly, we group $Y_{K+1},\ldots,Y_{2K}$ into another variable $Z_2$, and group $Y_{2K+1},\ldots,Y_{3K}$ into another variable $Z_3$. 
Then given the latent $\A$, the conditional probability table of $Z_1$, $Z_2$, $Z_3$ each has size $V^K \times H^K$; denote these three matrices by ${\bo\Psi}_m$ for $m=1,2,3$. According to the structure of the first $3K$ rows of the matrix $\G$, we can write ${\bo\Psi}_1$, ${\bo\Psi}_2$, and ${\bo\Psi}_3$ as follows:
$$
{\bo\Psi}_1
=\bigotimes_{j=1}^K \MP(Y_j\mid A_j),
\quad
{\bo\Psi}_2
=\bigotimes_{j=K+1}^{2K} \MP(Y_j\mid A_{j-K}),
\quad
{\bo\Psi}_3
=\bigotimes_{j=2K+1}^{3K} \MP(Y_j\mid A_{j-2K}).
$$
By our Assumption 1, each of the above three matrices is a Kronecker product of $K$ full column-rank matrices and hence also has full rank $H^K.$

Next, we further group the variable $Z_3$ and all the remaining variables $Y_{3K+1}, \ldots, Y_{J}$ (if they exist) into another categorical variable $Z_4$ with $V^{J-2K}$ categories. Denote the conditional probability table of $Z_4$ given $\A$ by ${\bo\Psi}_4$, which has size $V^{J-2K} \times H^K$. By definition, we have
$$
{\bo\Psi}_4 = {\bo\Psi}_3 \bigodot 
\underbrace{\MP(Y_{3K+1}\mid\A) \bigodot \MP(Y_{3K+2}\mid\A) \cdots \bigodot \MP(Y_{J}\mid\A)}_{J-3K\text{ matrices}}.
$$
Since every matrix in the above Khatri-Rao product is a conditional probability table with each column summing to one, the  ${\bo\Psi}_3$ can be obtained by summing appropriate rows of ${\bo\Psi}_4$. This indicates that the column rank of ${\bo\Psi}_4$ will not be smaller than that of ${\bo\Psi}_3$, so ${\bo\Psi}_4$ also has full rank $H^K$. 
Note that for alternative parameters $\ov\btheta$ and $\ov\nnu$, we have the following equations
\begin{align*}
    \Big({\bo\Psi}_1 \bigodot {\bo\Psi}_2 \bigodot {\bo\Psi}_4 \Big) \cdot \nnu = 
    \Big(\ov{\bo\Psi}_1 \bigodot \ov{\bo\Psi}_2 \bigodot \ov{\bo\Psi}_4 \Big) \cdot \ov\nnu 
\end{align*}
Now we invoke Kruskal's theorem \citep{kruskal1977three} as follows on the uniqueness of three-way tensor decompositions.
% \begin{lemma}[Kruskal's Theorem adapted from \cite{kruskal1977three}]\label{lem-kruskal}
Let $\mathbf M_1, \mathbf M_2, \mathbf M_3$ be three matrices of size $a_m\times r$ for $m=1,2,3$,  and $\mathbf N_1, \mathbf N_2, \mathbf N_3$ be three matrices each with $r$ columns. Suppose $\bigodot_{m=1}^3 \mathbf M_m  \cdot \one = \bigodot_{m=1}^3 \mathbf N_m \cdot \one$. Denote by $\rank_{\text{Kr}}(\mathbf M)$ the Kruskal rank of a matrix $\mathbf M$, which is the maximum number $R$ such that every $R$ columns of $\mathbf M$ are linearly independent.
If $\rank_{\text{Kr}}(\mathbf M_1) + \rank_{\text{Kr}}(\mathbf  M_2) + \rank_{\text{Kr}}(\mathbf  M_3) \geq 2r + 2$, then Kruskal's theorem guarantees that there exists a permutation matrix $\mathbf P$ and three invertible diagonal matrices $\mathbf D_m$ with $\mathbf D_1 \mathbf D_2 \mathbf D_3 =\mathbf I_r$ and $\mathbf N_m =\mathbf  M_m \mathbf D_m \mathbf P$ for each $m=1,2,3$.
% \end{lemma}

Based on Kruskal's theorem stated above, we obtain that $\bo\Psi_m=\ov{\bo\Psi}_m$ for $m=1,2,4$ and $\nnu=\ov\nnu$ up to a column latent class permutation. 
Finally, note that all entries of $\btheta$ can be directly read off from $\ov{\bo\Psi}_m$ for $m=1,2,3$ by definition. This implies the $\ov\btheta$ must also equal the $\btheta$ up to a label permutation. The proof of part (b) of the proposition is complete. 
\qed

\end{appendix}

\end{document}